\documentclass[final,leqno,onefignum,onetabnum]{siamltex1213}

\usepackage{graphicx} 
\usepackage{booktabs} 
\usepackage{amssymb}
\usepackage{pifont}
\usepackage{pgf,pgfarrows,pgfnodes,pgfautomata,pgfheaps,pgfshade}
\usepackage{tikz}
\usepackage{graphicx}
\usepackage{framed}
\usepackage{listings}
\usepackage{setspace}
\usepackage{mathtools}
\usepackage{enumitem}
\usepackage{subfigure}
\usepackage{caption}
\usepackage{multicol}

\begin{document}

\title{FEAST EIGENSOLVER FOR NON-HERMITIAN PROBLEMS} 

\author{James Kestyn\thanks{Department of Electrical and Computer Engineering, University of Massachusetts, Amherst, MA 01003, USA,
\email{kestyn@ecs.umass.edu}.}
\and
Eric Polizzi\thanks{Department of Electrical and Computer Engineering, Department of Mathematics and Statistics,
University of Massachusetts, Amherst, MA 01003 USA, \email{polizzi@ecs.umass.edu}.} 
\and
Ping Tak Peter Tang\thanks{Intel Corporation, Santa Clara CA 95054 USA, \email{peter.tang@intel.com}}}

\date{}

\maketitle
\slugger{sisc}{xxxx}{xx}{x}{x--x}


\begin{abstract}
A detailed new upgrade of the FEAST eigensolver targeting  
non-Hermitian eigenvalue problems is presented and thoroughly discussed.
It aims at broadening the class of eigenproblems that can be addressed within the 
framework of the FEAST algorithm.
The algorithm is ideally suited for computing
selected interior eigenvalues and their associated right/left bi-orthogonal eigenvectors, 
located within a subset of the complex plane. 
It combines subspace iteration with efficient contour integration 
techniques that approximate the left and right spectral projectors. 
We discuss the various algorithmic choices that have been 
made to improve the stability and usability of the new non-Hermitian eigensolver.
The latter
 retains the convergence property and multi-level parallelism of 
Hermitian FEAST, making it a valuable new software tool for the scientific community.
\end{abstract}



\pagestyle{myheadings}
\thispagestyle{plain}
\markboth{J. KESTYN, E. POLIZZI, AND P. TANG}{FEAST EIGENSOLVER FOR NON-HERMITIAN PROBLEMS}

\section{Introduction}

The generalized eigenvalue problem ${ AX=BX \Lambda}$  with $A$ and $B$ 
square matrices and $\Lambda$ diagonal,  is a central topic  
in numerical linear algebra and arises from a 
 broad and diverse set of disciplines in mathematics, science and engineering 
(the problem is said ``standard'' if $B \equiv I$ or ``generalized'' otherwise).
Solving the interior eigenvalue problem consists of determining
nontrivial solutions $\{\lambda_i, x_i\}$ (i.e. eigenpairs with $x_i=Xe_i$ and $\lambda_i=\Lambda_{i,i}$) 
located anywhere inside the spectrum.
Most common numerical applications lead to symmetric eigenvalue problems where $A$ is real symmetric or complex
Hermitian, $B$ is symmetric or Hermitian positive definite (hpd), and all the obtained eigenvalues $\lambda_i$ are real.
Non-symmetric and  non-Hermitian eigenvalue problem (including the case where A is complex symmetric) can also be encountered in a 
variety of situations resulting in complex values for $\lambda_i$. In this case $x_i$ is called the right eigenvector associated 
with $\lambda_i$, while one can also define a left eigenvector $\widehat{x}_i=\widehat{X}e_i$ solution of $\widehat{X}^HA=\Lambda\widehat{X}^HB$
(i.e.  $A^H\widehat{X}=B^H\widehat{X}\Lambda^*$).
Although many software packages are available for symmetric (or Hermitian) matrices (see e.g. \cite{hernandez2005survey,Arpack98,Sakurai2003,hernandez2005slepc,marek2014elpa,baker2009anasazi,stathopoulos2010primme,knyazev07blopex}), relatively few algorithms and software can handle 
the non-Hermitian problem  \cite{lehoucq1996evaluation,bai1997algorithm,Arpack98,bai1999able,garbow1978algorithm}. 
The FEAST eigensolver \cite{Polizzi09,FEASTsolver}, proven to be a robust and efficient tool for computing 
the partial eigenspectrum of Hermitian system matrices \cite{TangPolizzi14}, 
can also be generalized and applied to arbitrary non-Hermitian systems 
\cite{Laux2012,yin2014feast,TangKestynPolizzi14}.

FEAST is a subspace iteration method that uses 
the Rayleigh-Ritz projection and an approximate spectral projector as a filter \cite{TangPolizzi14}. 
Given a Hermitian generalized eigenvalue problem ${AX = BX \Lambda}$ of size $n$, the algorithm in Figure \ref{alg:feast} outlines the main
steps of a generic Rayleigh-Ritz subspace iteration procedure for computing $m$ eigenpairs.
\begin{figure}[htbp]
\fbox{
\begin{minipage}[t]{0.94\textwidth}
\noindent
 0. \ Start: Select random subspace $Y_{m_0} \equiv \{  y_1, y_2,\dots, y_{m_0}\}_{n\times m_0}$ \ ($n>>m_0\geq m$) \\  
\ 1. \ Repeat until convergence\\
\ 2. $\quad$   Compute  $ Q_{m_0} = \rho(B^{-1}A) Y_{m_0} $\\
\ 3. $\quad$   Orthogonalize ${ Q_{m_0}}$ \\ 
\ 4. $\quad$   Compute ${A_{Q}}=Q_{m_0}^H A Q_{m_0}$  and ${B_Q}=Q_{m_0}^HBQ_{m_0}$  \\
\ 5. $\quad$  Solve 
${{A_Q}{W}}={{B_Q}{W}{\Lambda_Q}}$ with ${W}^H{B_Q}{W}={I}_{m_0\times m_0}$\\ 
\ 6. $\quad$  Compute $ Y_{m_0}=Q_{m_0} W$ \\
\ 7. $\quad$ Check convergence of $Y_{m_0}$ and ${\Lambda_Q}_{m_0}$ for the $m$ 
wanted eigenvalues\\
\ 8. \ End 
\end{minipage}
}
\caption{Subspace iteration method with Rayleigh-Ritz projection
\label{alg:feast}} 
\end{figure} 
At convergence, the algorithm yields the $B$-orthonormal eigensubspace $Y_{m}\equiv X_{m}=\{  x_1, x_2,\dots, x_{m}\}_{n\times m}$ and 
associated eigenvalues  ${\Lambda_Q}_m\equiv{\Lambda}_m$. 
Taking  $\rho(B^{-1}A)=B^{-1}A$, yields the bare-bone subspace iteration (generalization of the power method)
which converges towards the $m$ dominant eigenvectors with 
the linear rate $|\lambda_{m_0+1}/\lambda_i|_{i=1,\dots,m}$ \cite{Saad89,Saad92,parlett1980symmetric}.
This standard approach is never used in practice.
Instead, it is combined with  filtering using the function $\rho$ which aims at improving the convergence rate  
(i.e. $|\rho(\lambda_{m_0+1})/\rho(\lambda_i)|_{i=1,\dots,m}$) 
by increasing the gap
between  wanted and  unwanted eigenvalues. 
The filtering function can also be expressed using the spectral decomposition of the Hermitian problem 
while considering the entire B-orthonormal eigensubspace i.e. $X^HBX=I$. 
\begin{eqnarray}
\rho(B^{-1}A) = X \rho(\Lambda) X ^{-1}\equiv X \rho(\Lambda) X ^{H}B.
\label{eq:spectral}
\end{eqnarray}
An ideal filter for the interior eigenvalue problem which maps all $m$ wanted eigenvalues 
to one and all unwanted ones to zero,
 can be derived from the
Cauchy (or Dunford)  integral formula:
\begin{eqnarray}
\rho(\lambda) = \frac{1}{2 \pi \imath} \oint_{\cal C} dz {(z - \lambda)^{-1}},
\label{eq:cauchy}
\end{eqnarray}
where the wanted eigenvalues are located inside a complex contour $\cal C$. 
The filter then becomes a spectral projector, with $\rho(B^{-1}A)=X_{m} X_{m} ^HB$, for the eigenvector subspace $X_m$ 
(i.e. $\rho(B^{-1}A)X_{m}=X_{m}$)  and can be written as: 
\begin{eqnarray}
\rho(B^{-1}A) = \frac{1}{2 \pi \imath} \oint_{\cal C}  dz (zB-A)^{-1}B.
\label{eq:density}
\end{eqnarray}
The FEAST method proposed in \cite{TangPolizzi14,Polizzi09},
uses a numerical quadrature to approximately compute the action of this filter onto a set of $m_0$ vectors along
the subspace iterations. The resulting rational function  $\rho_a$ that approximates the filter (\ref{eq:cauchy}) is given by:
\begin{eqnarray}
\rho_a(z) = \sum_{j=1}^{n_e} \frac{\omega_{j}}{z_j - z},
\label{eq:rational}
\end{eqnarray}
where $\{z_j,\omega_j\}_{1\leq j \leq n_e}$ are the nodes and related weights of the quadrature.
We obtain for the subspace $Q_{m_0}$ in step 2 of the algorithm in Figure \ref{alg:feast}: 
\begin{eqnarray}\label{eq:q0}
Q_{m_0}=\rho_a(B^{-1}A)Y_{m_0}= \sum_{j=1}^{n_e} {\omega_j}{(z_j B-A)}^{-1}BY_{m_0}\equiv X \rho_a(\Lambda) X ^{H}BY_{m_0}.
\end{eqnarray}
In practice, $Q_{m_0}$ can be computed by solving  a small number of (independent) shifted linear systems over a complex contour.
\begin{equation} \label{eq:q}
Q_{m_0} = \sum_{j=1}^{n_e} \omega_j Q^{(j)}_{m_0}, \ \quad \mbox{with $Q^{(j)}_{m_0}$ 
 solution of } \quad (z_j B-A)Q^{(j)}_{m_0}= BY_{m_0}
\end{equation}

The original FEAST paper \cite{Polizzi09}
demonstrated the effectiveness of the approach without analysis of convergence
or numerical issues. 
A detailed  numerical analysis on FEAST was completed 
recently in \cite{TangPolizzi14},
placing the algorithm on a more solid theoretical foundation.
In particular,  a relatively small number of quadrature nodes (using Gauss, Trapezoidal or Zolotarev \cite{guettel2014zolotarev} rules) on 
a circular contour suffices to produce a rapid decay of the function $\rho_a$
 from $\approx 1$ within the search contour to $\approx 0$ outside. 
In comparison with more standard polynomial filtering \cite{sorensen1992implicit,Saad92},
 the rational filter (\ref{eq:rational}) can lead to a very fast convergence of the subspace iteration procedure.
In addition, all the $m$ desired eigenvalues 
are expected to converge at the same rate (since $\rho_a(\lambda_i)\simeq 1$ if $\lambda_i$ is located within the search interval).
The convergence rate of FEAST 
does not only depend upon the decay properties of the rational function $\rho_a$, but 
also on the size of the search subspace $m_0$ which must not be chosen smaller than the number of eigenvalues inside the
search contour (i.e. $m_0\geq m$). 
Users of the FEAST eigensolver are responsible for specifying an interval to search for the eigenvalues and a subspace size $m_0$ 
that overestimate the number of the wanted eigenvalues.  Once these conditions are satisfied,
FEAST offers the following  set of appealing features: 
\begin{itemize}
\itemsep 1pt
\parskip 1pt
\item[(i)] high robustness with well-defined convergence rate  $|\rho_a(\lambda_{m_0+1})/\rho_a(\lambda_i)|_{i=1,\dots,m}$;
\item[(ii)] all multiplicities naturally captured; 
\item[(iii)] no explicit orthogonalization procedure on long vectors required in practice 
(i.e., step-3 in Figure \ref{alg:feast} is unnecessary as long as $B_Q$ is positive definite). 
We note in  (\ref{eq:q0}) that $Q_{m_0}$ is naturally spanned by the eigenvector subspace;
\item[(iv)] reusable subspace capable to generate suitable initial guess when solving a 
series of eigenvalue problems; 
\item[(v)] can exploit natural parallelism at
three different levels: search intervals can be treated separately (no overlap) while maintaining orthogonality - 
linear systems can be solved independently across the quadrature nodes of the complex contour - each complex linear system with $m_0$ 
multiple right-hand-sides can be solved in parallel. Consequently, within a parallel environment, the algorithm
complexity depends on solving a single linear system using a direct or an iterative method. 
\end{itemize}

By allowing the search contour to be placed at arbitrary locations in the complex plane,
the FEAST algorithm can be naturally extended to non-Hermitian problems which produce complex eigenvalues.
The algorithm retains most of the properties of Hermitian FEAST including 
the multi-level parallelism.
We note, however, a few theoretical and practical difficulties arising which
distinguish the non-Hermitian eigenvalue problems from Hermitian ones, including: (i) the treatment of defective systems using the
 Schur or Jordan forms; (ii) the notion of bi-orthogonality for dual right and left eigenvector subspaces; 
  (iii) the case of ill-conditioned eigenvalue problems that produce sensitive eigenvalues in finite precision arithmetic;
 (iv) or the shift-invert strategy that may give rise to ill-conditioned linear systems  
(e.g.  if a FEAST quadrature pole  lies near a complex eigenvalue).

 The key point at which the non-Hermitian FEAST algorithm differs from the Hermitian one is the use of dual subspaces. 
Since the left and right  eigenvectors do not necessarily lie in the same subspace, two separate projectors must then 
be calculated in order to  recover both sets of vectors. A single sided algorithm where only the right subspace 
is used to project is also  possible \cite{TangKestynPolizzi14}, but will not return a $B$-bi-orthogonal subspace of left and right eigenvectors, which can be of interest 
for many applications. In the following, all quantities associated
 with the left eigenvectors will be written with a '$\widehat{~~}$' symbol (e.g. $\widehat{X}$, $\widehat{Y}$ and $\widehat{Q}$). 
The non-Hermitian algorithm is similar to its Hermitian counterpart and follows the same steps outlined in Figure \ref{alg:feast}. 
 A comparison between the main numerical operations for the two algorithms is briefly outlined in 
Figure \ref{fig:2algo}. 
\begin{figure}[htbp]
\fbox{
\begin{footnotesize}
\begin{minipage}[t]{.47 \textwidth}
\underline{{\bf Hermitian FEAST}} \\

{\bf Solving:} $AX_m=BX_m\Lambda_m$ \\ 
 \hspace*{12.0mm}      $[{\Lambda_m}]_{ii}\in[\lambda_{min},\lambda_{max}]$ \\

{\bf Inputs:} {$A=A^H$, $B$ hpd; $m_0\geq m$; \\ 
\hspace*{12.0mm}$\{ z_1,\dots,z_{n_e} \}$, $\{ w_1,\dots,w_{n_e} \}$} \\

$Y_{m_0}$ $\leftarrow$ $m_0$ initial vectors \\
{\bf repeat}\\
 \hspace*{3mm} $Q_{m_0}=0$ \\
 \hspace*{3mm} {\bf for $j=1,n_e$} \\
 \hspace*{5mm}      $Q^{(j)}_{m_0} \leftarrow (z_j B - A)^{-1} BY_{m_0}$;  \\ 
 \hspace*{5mm}        $Q_{m_0} \leftarrow Q_{m_0} + \omega_j Q^{(j)}_{m_0} $ \\ 
 \hspace*{3mm} {\bf end} \\
 \hspace*{3mm} ${B_Q} \leftarrow Q^H_{m_0} B Q_{m_0}$ \\
 \hspace*{3mm} Check $B_Q$ hpd \ \ (resizing step)\\
 \hspace*{3mm} ${A_Q}\leftarrow Q^H_{m_0} A Q_{m_0}$ \\
 \hspace*{3mm} Solve ${A_Q}W = {B_Q} W {\Lambda_Q}; \ W^HB_QW=I$ \\
 \hspace*{3mm} $Y_{m_0} \leftarrow Q_{m_0} W$ \\
{\bf until} {\it Convergence of} $Y_m$, $\Lambda_{Q_m}$\\
\hspace*{8.0mm} {\it with}   $[{\Lambda_{Q_m}}]_{ii}\in[\lambda_{min},\lambda_{max}]$ \\

{\bf Output:} {$X_m\equiv Y_m$ ($X_m^HBX_m=I_m$); \\ 
\hspace*{12.5mm}  $\Lambda_m\equiv{\Lambda_Q}_m$}
\end{minipage}%
\hfill
\begin{minipage}[t]{.485\textwidth}
\underline{\bf Non-Hermitian FEAST}\\

{\bf Solving:} $AX_m=BX_m\Lambda_m ~~~~~~[{\Lambda_m}]_{ii}\in \cal C$ \\  $~~~~~~~~~~~~~~~A^H\widehat{X}_m=B^H\widehat{X}_m\Lambda_m^*$ \\

{\bf Inputs:}  {$A$ and $B$ general; $m_0\geq m$; \\
\hspace*{11.5mm} $\{ z_1,\dots,z_{n_e} \}$, $\{ w_1,\dots,w_{n_e} \}$}\\

 ${Y_{m_0}, \widehat{Y}_{m_0}}$ $\leftarrow$ $m_0$ initial vectors; \\
{\bf repeat}\\
 \hspace*{3mm} $Q_{m_0}=\widehat{Q}_{m_0}=0$  \\
 \hspace*{3mm} {\bf for $j=1,n_e$} \\
 \hspace*{5mm} $Q^{(j)}_{m_0} \leftarrow (z_j B - A)^{-1} BY_{m_0}$;  \\ 
 \hspace*{5mm} $\widehat{Q}^{(j)}_{m_0} \leftarrow (z_j^* B^H - A^H)^{-1} B^H \widehat{Y}_{m_0}$ \\
 \hspace*{5mm}  $Q_{m_0} \leftarrow Q_{m_0} + \omega_j Q^{(j)}_{m_0} $ \\
 \hspace*{5mm} $\widehat{Q}_{m_0} \leftarrow \widehat{Q}_{m_0} + \omega_j^* \widehat{Q}^{(j)}_{m_0} $\\
 \hspace*{3mm} {\bf end} \\
 \hspace*{3mm} ${B_Q} \leftarrow \widehat{Q}^H_{m_0} B Q_{m_0}$ \\
 \hspace*{3mm} Check $B_Q$ non-singular \ \ (resizing step)\\
 \hspace*{3mm} ${A_Q} \leftarrow \widehat{Q}^H_{m_0} A Q_{m_0}$ \\
 \hspace*{3mm} {Solve} ${A_Q}W={B_Q}W{\Lambda_Q} $ and \\
 \hspace*{3mm}  \hspace*{7mm} ${A_Q}^H\widehat{W}={B_Q}^H \widehat{W}\Lambda_Q^*; \ \widehat{W}^HB_QW=I$ \\
 \hspace*{3mm} $Y_{m_0} \leftarrow Q_{m_0} W, \quad \widehat{Y}_{m_0}\leftarrow\widehat{Q}_{m_0}\widehat{W}$ \\ 
{\bf until} {\it Convergence of}  $Y_m, \widehat{Y}_m, {\Lambda_Q}_m$ \\
\hspace*{8.0mm}{\it with } $[{\Lambda_{Q_m}}]_{ii}\in \cal C$ \\


{\bf Output:} {$X_m\equiv Y_m$; \\ \hspace*{12.25mm} $\widehat{X}_m \equiv\widehat{Y}_m$ ($\widehat{X}_m^HBX_m=I_m$); \\ 
\hspace*{12.75mm} $\Lambda_m\equiv{\Lambda_Q}_m$}
\end{minipage}%
\end{footnotesize}
}
\caption{\label{fig:2algo} Brief outlook and comparison between the main numerical operations for the FEAST algorithms 
applied to the Hermitian and non-Hermitian problems.}
\end{figure}
The rest of the article aims at providing all the details of the non-Hermitian FEAST algorithm and its practical implementation.
Section~\ref{sec:2} presents multiples theoretical and practical algorithmic considerations, outlines the differences with the Hermitian 
FEAST algorithm, and ends with a 
complete description  of the non-Hermitian algorithm with discussions on limitations.
 Section~\ref{sec:3} briefly outlines some features of the
 new FEAST eigensolver version 3.0,  from which the proposed changes here take effect.
We conclude by presenting some numerical experiments in Section~\ref{sec:4}. 


\section{Theoretical and Practical Considerations}\label{sec:2}

\subsection{Defining a search contour}\label{sec:contour0}

\begin{figure}[htbp]
 \center{\includegraphics[width=1.0\textwidth]
        {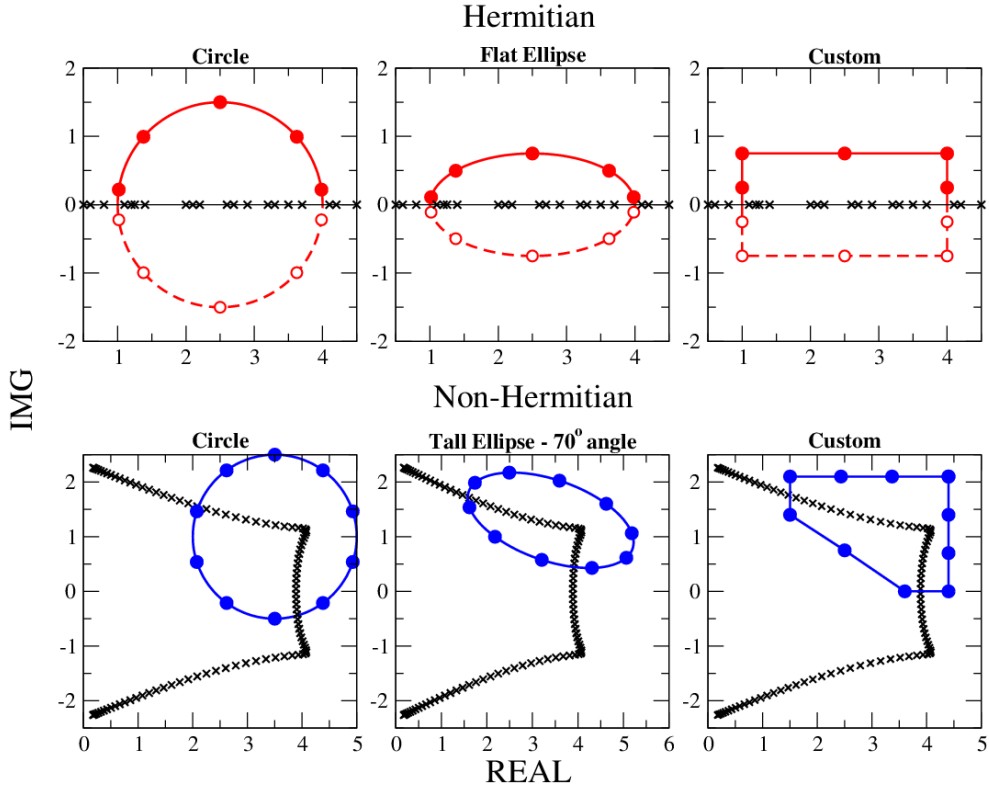}}
        \caption{\label{fig_contour} Various search contour examples for the Hermitian and the non-Hermitian FEAST algorithms.
Both algorithms feature standard ellipsoid contour options and the possibility to define custom arbitrary shapes. 
In the Hermitian case, the contour is symmetric with the real axis and only the nodes in the upper-half may be generated.
In the non-Hermitian case, a full contour is needed to enclose the wanted complex eigenvalues. 
}
\end{figure}

A key feature of FEAST is the ability to calculate a subset of eigenvalues that exist within some interval. 
Figure \ref{fig_contour} summarizes the different search contour options possible 
for both the Hermitian and non-Hermitian FEAST algorithms.

For the Hermitian case, the user must then specify a 1-dimensional real-valued search interval $[\lambda _{min}, \lambda _{max}]$. 
These two points are used to define a circular or ellipsoid contour $\cal C$ centered on the real axis, and along
 which the complex integration nodes are generated. 
The choice of a particular quadrature rule 
will lead to a different set of relative positions for the nodes and associated quadrature weights i.e.  $\{ z_j,\omega_j\}$.
Since the eigenvalues are real, it is convenient to select a symmetric contour 
with the real axis (i.e. ${\cal C}= {\cal C}^*$) since it only requires one to operate the quadrature
on the half-contour (e.g. upper half). 

With a non-Hermitian problem, it is necessary to specify  a 2-dimensional search contour that surrounds
the wanted complex eigenvalues. 
Circular or ellipsoid contours can also be used and they can be generated using standard options included into FEAST v3.0. 
These are defined by a complex midpoint $\lambda_{mid}$ and a radius $r$ for a circle (for an ellipse the ratio 
between the horizontal axis and vertical axis diameter can also be specified, as well as an angle of rotation). However,
in some applications where the eigenvalues of interest belong to a particular subset in the complex plane,  
more flexibility for selecting a search contour with arbitrary shape could be needed.
This option also  lends itself to parallelism, where a large number of eigenvalues can be calculated by partitioning the complex 
plane into multiple contours (see Section \ref{sec:4}).
Consequently, a ``Custom Contour'' feature is also supported in FEAST v3.0 that  allows to account for 
arbitrary 
 quadrature nodes and weights.

\subsection{Right/Left Spectral Projectors and Dual Subspaces}

The filtering function 
can be applied to any similarity transformation of the pencil, the most general of which is the Jordan Normal Form. 
\begin{equation}
\rho(B^{-1}A) = X \rho(J) X^{-1}.
\label{eq:jordan0}
\end{equation}
When applied to each Jordan block $J_k$, the expression of the operator becomes \cite{higham2008functions}:
\begin{equation}
\rho(J_k) = \rho \left(
\begin{bmatrix}
    \lambda _{k} & 1        & \dots  & 0 \\
    0     & \lambda_{k}    & \ddots  & \vdots \\
    \vdots & \ddots  & \ddots & 1 \\
    0      & \dots   & 0      & \lambda_{k}
\end{bmatrix} 
\right) =
\left[
\arraycolsep=1.4pt\def\arraystretch{2.2}
\begin{array}{cccc}
    \rho(\lambda _{k}) & \dfrac{ \rho^{\prime} (\lambda _{k}) }{0!}       & \dots  & \dfrac{ \rho^{(m)} (\lambda _{k}) }{(m-1)!} \\
    0     & \rho( \lambda_{k} )   & \ddots  & \vdots \\
    \vdots & \ddots  & \ddots &  \dfrac{ \rho^{\prime} (\lambda _{k}) }{0!}  \\
    0      & \dots   & 0      & \rho(\lambda_{k})
\end{array}
\right]
\label{eq:jordan1}
\end{equation}
Using the Cauchy integral formula (\ref{eq:cauchy}), the diagonal elements of $J_k$ take the values one or zero, while
all derivatives (i.e. off-diagonal elements) are zero. In practice,  this may not be guarantee with FEAST as the filter
 is approximated by the rational function (\ref{eq:rational}). 
A generalization of the algorithm for addressing the defective systems would require further studies, and our current
FEAST non-Hermitian algorithm assumes that the Jordan form reduces to an eigenvalue decomposition.
 Consequently, we consider: 
\begin{eqnarray}
\rho(B^{-1}A) = X \rho(\Lambda) X ^{-1}\equiv X \rho(\Lambda) \widehat{X} ^{H}B,
\label{eq:spectral2}
\end{eqnarray}
where the left and right eigensubspaces satisfy the $B$-bi-orthonormal relationship i.e. 
$\widehat{X}^H B X = I$. For the case of the Hermitian problem, we note that $\widehat{X}={X}$ and  
the relation (\ref{eq:spectral}) can then be recovered. It is also important to mention the particular case
of complex symmetric systems (i.e. $A=A^T$ and $B=B^T$) which leads to  $\widehat{X}={X}^*$.
In general, however, the left and right vectors are not straightforwardly related and they must be calculated explicitly.

From (\ref{eq:cauchy}), (\ref{eq:density}), 
and (\ref{eq:spectral2}), one can define the right spectral projector $\rho(B^{-1}A)$ for  
the right eigenvector subspace $X_m$ (i.e. $\rho(B^{-1}A)X_{m}=X_{m}$) as follow:
\begin{eqnarray}
\rho(B^{-1}A) = \frac{1}{2 \pi \imath} \oint_{\cal C}  dz (zB-A)^{-1}B \equiv X_{m}\widehat{X}_{m}^HB.
\label{eq:right}
\end{eqnarray}
For the treatment of the left eigenvector subspace solution of $\widehat{X}^HA=\Lambda\widehat{X}^HB$, 
it is first convenient to define the following eigenvalue decomposition: 
\begin{eqnarray}
\rho(AB^{-1}) = \widehat{X}^{-H} \rho(\Lambda) \widehat{X}^H\equiv B X \rho(\Lambda) \widehat{X}^{H}.
\label{eq:spectral2hat}
\end{eqnarray}
One can then construct the left spectral projector $\rho(AB^{-1})$ 
 (i.e. $\widehat{X}_m^H=\widehat{X}_m^H\rho(AB^{-1})$) as:
\begin{eqnarray}
{\rho}(AB^{-1}) = \frac{1}{2 \pi \imath} \oint_{\cal C}  dz B(zB-A)^{-1} \equiv B {X}_{m}\widehat{X}_{m}^H.
\label{eq:left}
\end{eqnarray}
%
In FEAST, the projectors 
are formulated using the rational function $\rho_a$ (\ref{eq:rational})  along with
the quadrature nodes and weights $\{z_j,\omega_j\}_{1\leq j \leq n_e}$ that approximate the contour integrations in 
(\ref{eq:right}) and (\ref{eq:left}).
The right and left subspaces $Q_{m_0}$ and $\widehat{Q}_{m_0}$ are then obtained by applying the right and left projectors 
onto a set of $m_0$ vectors i.e.
\begin{eqnarray}\label{eq:q00}
{Q}_{m_0}={\rho}_a(B^{-1}A){Y}_{m_0}= \sum_{j=1}^{n_e} {\omega_j}{(z_j B-A)}^{-1}B{Y}_{m_0}\equiv 
{X} \rho_a(\Lambda) \widehat{X}^{H}B{Y}_{m_0}.
\end{eqnarray}
and 
\begin{eqnarray}\label{eq:qhat0}
\widehat{Q}_{m_0}^H=\widehat{Y}_{m_0}^H{\rho}_a(AB^{-1})= \sum_{j=1}^{n_e} {\omega_j}\widehat{Y}^H_{m_0}B{(z_jB-A)}^{-1}\equiv 
\widehat{Y}_{m_0}^HB{X} \rho_a(\Lambda) \widehat{X}^{H}.
\end{eqnarray}
In practice, the calculation of both subspaces require solving a series of linear systems. For the right subspace, 
\begin{equation} 
\label{eq:qr}
Q_{m_0} = \sum_{j=1}^{n_e} \omega_j Q^{(j)}_{m_0}, \ \quad \mbox{with $Q^{(j)}_{m_0}$ 
 solution of } \quad (z_j B-A)Q^{(j)}_{m_0}= BY_{m_0},
\end{equation} 
which was already outlined  in (\ref{eq:q}), and for the left subspace:
\begin{equation} \label{eq:ql}
\widehat{Q}_{m_0} = \sum_{j=1}^{n_e} \omega_j^* \widehat{Q}^{(j)}_{m_0}, 
\ \quad \mbox{with $\widehat{Q}^{(j)}_{m_0}$ 
 solution of } \quad 
(z_jB-A)^H\widehat{Q}^{(j)}_{m_0}= B^H\widehat{Y}_{m_0}.
\end{equation} 
These numerical operations are also described in Figure \ref{fig:2algo}. 
As a result of (\ref{eq:q00}) and (\ref{eq:qhat0}), 
$Q_{m_0}$ (resp. $\widehat{Q}_{m_0}$) is formed by a linear combinations of the columns of $X_{m_0}$ 
(resp. $\widehat{X}_{m_0}$).  The Rayleigh-Ritz procedure should then involve
the reduced matrices 
${A_Q}$ and ${B_Q}$ formed by projecting on the right with a subspace containing the right eigenvectors 
$Q_{m_0}$, and projecting  on the left with a subspace containing the left eigenvectors $\widehat{Q}_{m_0}$. 
The resulting non-Hermitian reduced system can be solved using the QZ algorithm \cite{moler1973algorithm} in LAPACK \cite{anderson1999lapack} to yield
the right and left eigenvectors $W$ and $\widehat{W}$ defined in Figure \ref{fig:2algo}. 
The long right (resp. left) Ritz vectors can then be recovered as ${Y}_{m_0}={Q}_{m_0}W$ (resp. $\widehat{Y}_{m_0}=\widehat{Q}_{m_0}\widehat{W}$), 
and used as initial guess subspaces for the next FEAST iterations until convergence.


\subsection{Discussions on Convergence}\label{sec:convergence}

In our implementation of FEAST, the criteria of convergence is satisfied if
the norm of the relative residual associated with the eigenpairs $(x_i, \lambda_i)$ and  $(\widehat{x}_i, \lambda_i^*)$, 
is found below an arbitrary threshold $\epsilon$ i.e.
\begin{equation}
\label{eq:res}
res_i=\max
\left\{\frac{||Ax_i - \lambda_i Bx_i ||_1}{||\alpha Bx_i ||_1},\frac{||A^H\widehat{x}_i - \lambda_i^* B^H\widehat{x}_i ||_1 }
{||\alpha B^H\widehat{x}_i ||_1}\right\}<\epsilon,  
\end{equation}
where the value of
$\epsilon$ can be chosen typically equal to $10^{-13}$ if high accuracy is needed using double precision arithmetic.
The parameter $\alpha$ is relative to the eigenvalue range in the search contour. The latter  is defined differently for the Hermitian 
and non-Hermitian cases (as discussed in Section \ref{sec:contour0}), and  a non-zero value for $\alpha$ can be chosen as 
$\alpha = \max(|\lambda_{min}|, |\lambda_{max}|)$ for the Hermitian case and $\alpha = (|\lambda_{mid}| + r)$ for the non-Hermitian case.

As discussed in the introduction section, the right/left eigenvectors associated with $\lambda_i$ with $i=1,\dots,m$ 
(and hence all associated residuals $res_i$) are expected to converge linearly along the FEAST subspace iterations at the rate:  
$|\rho_a(\lambda_{m_0+1})/\rho_a(\lambda_i)|$ for $i=1,\dots,m$. 
The convergence depends then on both the subspace size $m_0$ ($m_0\geq m$) and the accuracy of the rational filter $\rho_a$
(\ref{eq:rational}) that should ideally provide values very close to unity for eigenvalues on the interior 
of the search contour and zero elsewhere.
Although quite effective, the Gauss-quadrature approach along a circular contour that was proposed in the original FEAST 
article \cite{Polizzi09}, is clearly not the only possible choice for optimizing the convegence ratio. 
Three other options have already been considered for the Hermitian problem including \cite{guettel2014zolotarev}: (i) the Trapezoidal rule; (ii) different contour
shapes beside a circle such as a finely tuned flat ellipse; (iii) a new approximation of the spectral projector
based on a Zolotarev approximant to the sign function which, after transformations, provides complex
poles on the unit circle \cite{zolotarev1877application,guettel2014zolotarev}. Both Gauss and Zolotarev are well-suited choices for the Hermitian problem since they
favor an accentuation of
the decay of $|\rho_a|$ at the boundaries of the interval along the real axis. The Trapezoidal rule, in turn,
leads to a more uniform decay for  $|\rho_a|$  in any directions of the complex plane \cite{TangKestynPolizzi14}, and it 
is also well-known for its exponential convergence property with the number of integration nodes $n_e$ \cite{trefethen2014exponentially}.
The Trapezoidal rule is then expected to provide more consistency
for capturing the complex eigenvalues of the non-Hermitian problem.

\begin{figure}[htbp]
 \center{\includegraphics[width=0.95\textwidth]
        {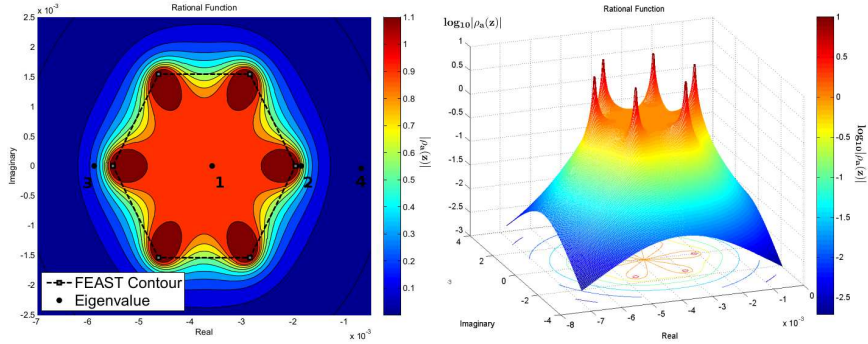} }
        \caption{\label{fig_ratf} Value of the rational function plotted as contour plot (left) and surface plot (right) using a
 hexagonal contour for the QC324 matrix. The  left plot includes the positions of the four closest eigenvalues. Only a single 
eigenvalue $\lambda_1$ is inside of the contour. More particularly, we note that
 $|\rho_a(\lambda _1)| =1.0000004$ , $|\rho_a(\lambda _2)| = 1.7272309$, $|\rho_a(\lambda _3)| = 0.4206553$, 
$|\rho_a(\lambda _4)| =3.6296209\times 10^{-2}$, and $|\rho_a(\lambda _5)| =6.9332547\times 10^{-3}$. 
The latter is associated with the eigenvalue $\lambda_5$ which cannot be seen in the Figure since it is out of range.
}
\end{figure}

\begin{figure}[htbp]
 \center{\includegraphics[width=0.9\textwidth]
        {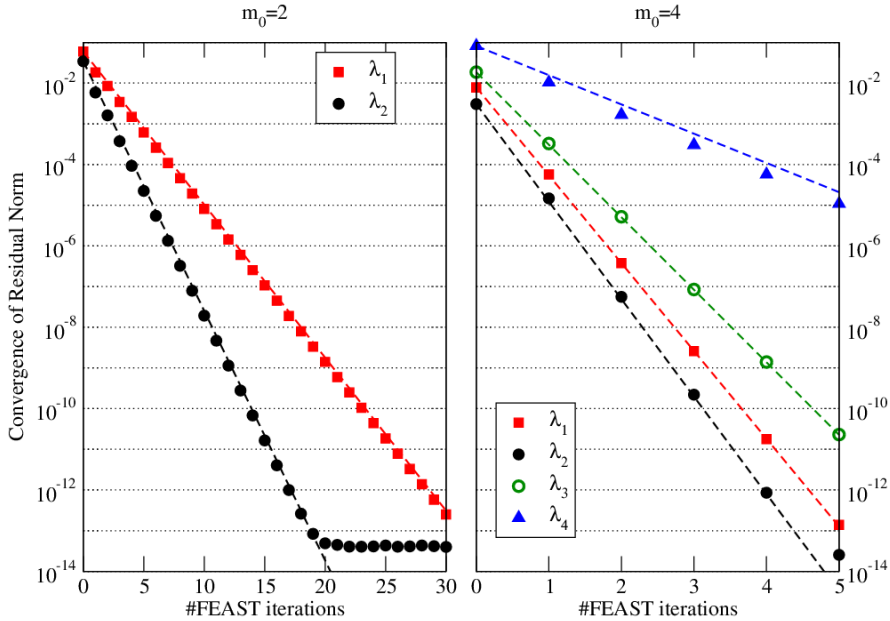}}
        \caption{\label{fig_convg} Convergence of the residual norms 
(\ref{eq:res}) associated with eigenvalues $\lambda_i$ in Figure \ref{fig_ratf}.
Two search subspace size are considered $m_0=2$ (left plot) and $m_0=4$ (right plot). The dashed lines represent
the theoretical linear convergence rate  $|\rho_a(\lambda_{m_0+1})/\rho_a(\lambda_i)|$ which is perfectly matched  
by the values returned by FEAST. We note that the convergence of the wanted eigenvalue $\lambda_1$ is
 is considerably slower using the smaller size subspace $m_0=2$ since the eigenvalue $\lambda_3$, that governs the convergence rate 
for this case, ends up being too close to the search contour.         }
\end{figure}

Similarly to the Hermitian problem, the non-Hermitian FEAST algorithm also requires as input
 a search subspace of size $m_0$ chosen not smaller than the number of eigenvalues $m$ within a given complex contour. 
If $m_0$ ($m_0\geq m$) is chosen too small, the ratio governing the convergence may come closer to one, leading to slow convergence.
Alternatively, if $m_0$ is chosen too large, it may result in an unnecessary high number of right-hand-sides
 when solving the shifted linear systems in (\ref{eq:qr}) and (\ref{eq:ql}).  
Two examples have been designed to illustrate the convergence rates dependence on $m_0$. 
These tests use the QC324 matrix from the NEP collection \cite{Bai96}. A contour has been created with a single eigenvalue $\lambda_1$ inside.
The contour and few closest eigenvalues can be seen in Figure \ref{fig_ratf}. The rational function $|\rho_a|$ which 
has been generated using a six-point Trapezoidal rule is also shown in the figure
 (as a contour plot on the left and a 3-D surface plot on the right).   
Figure \ref{fig_convg} shows the convergence of the relative residual norms for all 
the $m_0$ eigenvalues along the FEAST subspace iterations in the cases $m_0=2$ (left plot) and $m_0=4$ (right plot).
For $m_0=2$, $\lambda_3$ is the closest eigenvalue outside of the search subspace and controls
 the convergence rate. Since this eigenvalue is relatively close to the contour, 
FEAST exhibits slow convergence. The case $m_0=4$, in turn, leads to drastic improvement in the
  convergence rate which benefits from the small values of $\rho_a(\lambda_5)$.


A typical recommended choice for the search subspace size is $m_o=2m$. 
In practice, however, the exact number of eigenvalues $m$ is unknown beforehand and the user must make an educated guess.
Alternatively, $m$ can also be estimated using, for example,
the fast stochastic estimate procedure \cite{di2013efficient} that has been recently introduced in FEAST v3.0. 
It is important to note that in some situations 
slow convergence can result if the value of $m_0$ is only large enough to include the external eigenvalues bordering the contour.
This problem can arise when the eigenvalues of interest
 are near a continuum or cluster of eigenvalues. With many eigenvalues closely bordering the contour, it may not
 be possible to improve convergence by increasing the subspace size $m_0$. In this case, using additional integration 
nodes to increase the accuracy of $\rho_a$ may be necessary. A utility routine for calculating the rational 
function has also been included in FEAST v3.0 and can be used to investigate convergence for different 
contours and eigenvalue distributions. 

Finally, and in contrast to the Hermitian problem where the contour nodes can be placed away from the eigenvalues 
(i.e. far enough from the real axis), 
 a contour node could end up being located in the vicinity of a complex eigenvalue.
In this case the rational function could take on values larger than one, and it becomes then possible
 for an eigenvalue outside of the contour to converge at a faster rate 
 than the wanted eigenvalues inside. This is what is happening to $\lambda_2$ in Figures  \ref{fig_ratf} and \ref{fig_convg}.
If a contour node is located too close to an eigenvalue, however, 
it is likely to worsen the conditioning of the corresponding shifted linear system in (\ref{eq:qr}) and (\ref{eq:ql}), 
making then the problem  more challenging to solve using an iterative method.


\subsection{Reduced Contour Integration Cost}\label{sec:contour} 

Non-Hermitian matrices $A$ and $B$ come in three flavors: (i) complex general, (ii) real non-symmetric, 
and (iii) complex symmetric. The major computational task performed by FEAST is the numerical integration,
 where a set of linear systems must be solved  along a complex contour. 
In the  complex general case both $Q_{m_0}$ and $\widehat{Q}_{m_0}$ are 
computed explicitly
by solving the $2n_e$ (independent) linear systems defined in (\ref{eq:qr}) and (\ref{eq:ql}).
It is important to note that most modern numerical libraries that 
includes direct methods for solving linear systems, supply a ``transpose conjugate solve'' feature as well 
(i.e. a linear system $A^Hx=f$ can be solved using the factorization of $A$).  
Consequently, once the $(z_jB-A)$ matrices are factorized in (\ref{eq:qr}), the system solves in (\ref{eq:ql}) 
can be performed without re-factorizing the conjugate transpose of the matrices. Similarly using iterative methods, 
the conjugate transpose solve could be performed without factorizing twice the preconditioner.
If such option is available, the contour integration in  the most general case
should involve only $n_e$ (independent) factorizations 
and $2n_e$ (independent) solves with $m_0$ right hand sides. 
For the cases (ii) and (iii) above,   it is
possible to take advantage of some additional matrix properties that result in a reduced workload as discussed in following.
\\
\begin{description}
\itemsep 1pt
\parskip 1pt
\item[Complex symmetric]-
 For the complex symmetric case  (${A} = {A}^T$ and $B=B^T$), there exists
 a relationship between the left and right eigenvectors, which can be expressed as conjugate pairs i.e. $\widehat{X} = X^*$. This allows the left  
subspace $\widehat{Q}$ to be expressed in terms of the right $Q$ using the same simple relationship $\widehat{Q} = Q^*$.
Therefore $\widehat{Q}$ (\ref{eq:ql}) does not need to be calculated, and 
only the $n_e$ factorizations and $n_e$ solves in (\ref{eq:qr}) are then necessary. \\ 

\item[Real non-symmetric]-
In general the treatment of the real non-symmetric case  (${A} = {A}^*$ and $B=B^*$) is identical to the
complex non-symmetric one. However, there exists some savings for specific contours exhibiting symmetry across the real axis 
(i.e  $\cal C=\cal C^*$). 
For this particular case, each integration node $z_j$ with $j=1,\dots,n_e/2$ 
in the upper half of the complex plane has a conjugate pair $z_j^*$  
in the lower half. From the resulting following relationships:
$$(z_jB-A)^*=(z_j^*B-A) \quad \mbox{and} \quad (z_jB-A)^H=(z_j^*B-A)^T,$$
one can show that only the $n_e/2$ factorizations of $(z_jB-A)$ in the upper-half contour, along with
 $n_e$ total solves, are needed to obtain both $Q_{m_0}$ and $\widehat{Q}_{m_0}$ in (\ref{eq:qr}) and (\ref{eq:ql}).



\end{description}

The contour integration cost can then be reduced depending on the properties of the eigenvalue system, attributes
of the complex contour (e.g. if $\cal C=\cal C^*$), or the standard feature of transpose conjugate solve offered by most
 linear system solvers.
Table \ref{tab:solve} summarizes the number of factorizations and solves effectively needed to perform
the full contour integration using a total of $n_e$ nodes. 
The cost of the Hermitian FEAST algorithm is also provided for reference.

\begin{table}[htbp]
\begin{small}
\begin{center}
\begin{tabular}{llcc}
Family of eigenvalue problems  & ($A$,$B$) properties &\#Factorizations & \#Solves \\
\hline \\[-1.5ex]
Complex general & N/A & $n_e$ & $2n_e$ \\
Complex symmetric & $A=A^T$, $B=B^T$ &$n_e$ & $n_e$\\
Complex Hermitian with $\cal C=\cal C^*$ & $A=A^H$, $B$ hpd & $n_e/2$ & $n_e$ \\
Real non-symmetric & N/A & $n_e$ & $2n_e$ \\
Real non-symmetric with $\cal C=\cal C^*$ & N/A & $n_e/2$ & $n_e$ \\
Real symmetric with $\cal C=\cal C^*$ & $A=A^T$, $B$ spd & $n_e/2$ & $n_e/2$ \\
\hline
\end{tabular}
\caption{\label{tab:solve} Summary of the total number of factorizations and solves effectively needed by FEAST to perform
the full contour integration using a total of $n_e$ nodes. 
It is also assumed that the transpose conjugate solve feature is available for the system solver.
}
\end{center}
\end{small}
\end{table}


\subsection{Resizing the search subspace}\label{sec:resize}



The rank of the subspaces $Q_{m_0}$ (\ref{eq:q00}) and $\widehat{Q}_{m_0}$  (\ref{eq:qhat0})
is greater than or equal to the number of wanted eigenvalues $m$ ($m_0\geq m$) since 
the eigenpairs outside of the contour are also accounted for due to inaccuracies in the numerical integration. 
In turn, if $m_0$ is too  severely overestimated 
the rank may end up being less than the subspace size $m_0$ in finite precision arithmetic.
Consequently, 
 $m_0$ 
 must then be resized
to prevent the subspaces to become numerically rank deficient and the reduced matrix 
${B_Q}=\widehat{Q}_{m_0}^HBQ_{m_0}$   singular.
Otherwise, the QZ algorithm used in the computation of the reduced system can produce infinite eigenvalue solutions \cite{moler1973algorithm}.
 Re-injecting these solutions into the subspace iteration would cause problems for the algorithm. 
The upper bound for the choice of $m_0$ should be the largest value before the subspaces
become numerically rank deficient.
One possible way to determine 
this threshold value consists of
performing 
the spectral decomposition of ${B_Q}$ and analyzing its eigenvalues. It comes:
\begin{eqnarray}
\label{eq:bdiag}
{B_Q}= V \Gamma \widehat{V} ^H,
\end{eqnarray} 
where ${\Gamma}$  is the diagonal matrix for the eigenvalues $\{\gamma_i\}_{i=1,\dots,m_0}$, and   $V$ and $\widehat{V}$ are respectively
the corresponding left and right bi-orthonormal eigenvector subspaces (i.e. $\widehat{V}^HV=I_{m_0}$).
We note that for the Hermitian case where $B$ and hence ${B_Q}$ must be positive definite, 
this step is replaced by monitoring the failure of the Cholesky factorization of ${B_Q}$ that could return a negative pivot.
The position of the latter helped determining the threshold value for $m_0$ used to resize the subspace accordingly.
For the non-Hermitian problem, the matrix ${B_Q}$ is singular if there exists an eigenvalue equal to zero.
 In finite precision arithmetic, a zero eigenvalue must be characterized relatively i.e.
\begin{equation}
\label{eq:resize}
 |\gamma_i| < \eta * \max\left\{|\gamma _1|, ... ,|\gamma_{m_0}|\right\},
\end{equation} 
where $\eta$ is relative to the machine precision; e.g. $10^{-16}$ in double precision.
If an eigenvalue $\gamma _i$ is then different than the maximum eigenvalue by 16 orders-of-magnitude 
then it is out of range for the double precision arithmetic and is counted as a zero. 
The subspace $m_0$ is resized to $\widetilde{m}_0$ such that ${B_Q}$ has no eigenvalues satisfying (\ref{eq:resize}).
The spectral decomposition of ${B_Q}$ is computed at each FEAST iteration, and as it will be discussed in the next section,
 the resizing is performed in conjunction with a $B$-bi-orthonormalization for
the subspaces $Q_{\widetilde{m}_0}$ and $\widehat{Q}_{\widetilde{m}_0}$. 
The additional numerical cost of diagonalizing ${B_Q}$ is on the order of (but less expensive than)
 the cost associated with the diagonalization of the reduced generalized system.

As a side remark, it interesting to note that
using the expression (\ref{eq:q00}) and (\ref{eq:qhat0}),  $B_Q$ can also be written as:
\begin{eqnarray}
\label{eq:bdiag1}
{B_Q}= (\widehat{Y}_{m_0}^HBX) \rho_a^2(\Lambda) (\widehat{X}^HBY_{m_0}).
\end{eqnarray} 
Starting from the second FEAST iteration where the Ritz vectors $Y_{m_0}$ and $\widehat{Y}_{m_0}$  are not only span respectively by the  true 
eigenvector subspaces $X$ and $\widehat{X}$ but they also satisfy
 the property of $B$-bi-orthonormality (i.e. $\widehat{Y}_{m_0}^HBY_{m_0}=I$ since $\widehat{W}^HB_QW=I$ in Figure \ref{fig:2algo}), 
it is possible to directly identify (\ref{eq:bdiag1}) with (\ref{eq:bdiag}).
It comes that $V=\widehat{Y}_{m_0}^HBX$, $\widehat{V}^H=\widehat{X}^HBY_{m_0}$, and $\Gamma=\rho_a^2(\Lambda)$. The latter indicates
that the eigenvalues of $B_Q$ are related to the rational function $\rho_a$, and can then be used to 
estimate the convergence rate \cite{TangPolizzi14}.
In order for $|\rho_a|$ to satisfy (\ref{eq:resize}), however, $\eta$ should be replaced by $\sqrt{\eta}$. Consequently, the convergence
rate for the algorithm is here limited to $10^{-8}$ in double precision arithmetic 
(a similar argument could be made for the case of the 
Hermitian FEAST which relies on the Cholesky decomposition of the normal-type equation $B_Q$).
FEAST can then converge in a minimum of 2 iterations to machine precision
 at $\sim10^{-16}$ given a sufficiently large enough subspace size $m_0$ (whose value is also 
relative to the accuracy of $\rho_a$). If needed, it may be possible to obtain higher convergence rate (i.e one FEAST iteration)
using a direct robust QR factorization or singular value decomposition (SVD) of the subspaces
 $Q_{m_0}$ and $\widehat{Q}_{m_0}$. 

\subsection{B-bi-orthonormalization}\label{sec:biortho}

The intended result of FEAST is a set of $B$-bi-orthonormal vectors. 
However, the $B$-bi-orthogonality is not guaranteed after the contour integration due to numerical inaccuracies.
 This is especially pronounced in large problems which exhibit a continuum of eigenvalues bordering the search contour. 
The contour integration could potentially include a large number of mixed states from the continuum 
in the subspaces $Q_{m_0}$ (\ref{eq:q00}) 
and $\widehat{Q}_{m_0}$ (\ref{eq:qhat0}). 
In our numerical experiments, we have found that an explicit $B$-bi-orthonormalization of the FEAST subspaces
 $Q_{m_0}$ and $\widehat{Q}_{m_0}$ helps improving the stability of the algorithm.
Rather than performing a QR factorization or SVD of the subspaces,
we aim at taking advantage of the eigen-decomposition of ${B_Q}$ (\ref{eq:bdiag}) that is 
 already performed in FEAST as discussed in the previous section. 
From  (\ref{eq:bdiag}) and since $B_Q= \widehat{Q}^H_{m_0} B Q_{m_0}$, it comes:
\begin{eqnarray}
\Gamma = \widehat{V}^H {B_Q} V = (\widehat{V}^H \widehat{Q}^H_{m_0}) B (Q_{m_0} V)\equiv (\widehat{Q}_{m_0}\widehat{V})^H B (Q_{m_0} V).
\end{eqnarray} 
As a result,  $B$-bi-orthonormal subspaces $U_{m_0}$ and $\widehat{U}_{m_0}$ can be generated by updating the current subspaces 
$Q_{m_0}$ and $\widehat{Q}_{m_0}$ as follows:
\begin{eqnarray}
\label{eq:u}
U_{m_0} = Q_{m_0} V \Gamma^{-1/2}, \hspace*{10mm} \widehat{U}_{m_0} = \widehat{Q}_{m_0} \widehat{V}{\Gamma}^{-H/2}.
\end{eqnarray}
As discussed in the previous section, the subspace size $m_0$ may have already been reduced to $\widetilde{m}_0$ at this stage by
allowing the eigenvectors in  $V$ and $\widehat{V}$, 
corresponding to the zero eigenvalues in $\Gamma$, to be removed from the subspace. 
In practice, a subset of $V$ and $\widehat{V}$ composed of $\widetilde{m}_0$ columns vectors 
 can be easily extracted if the eigenpairs $\{\gamma_i,v_i\equiv Ve_i, \widehat{v}_i\equiv \widehat{V}e_i\}_{i=1,\dots,m_0}$ 
are first sorted by decreasing values of $|\gamma_i|$. Denoting $V_{m_0\times \widetilde{m}_0}$
and $\widehat{V}_{m_0\times \widetilde{m}_0}$ the subsets of the new $V$ and $\widehat{V}$ subspaces restricted to their first 
 $\widetilde{m}_0$ columns, and $\Gamma_{ \widetilde{m}_0\times \widetilde{m}_0}$ the matrix of the first  $\widetilde{m}_0$ sorted eigenvalues,
 (\ref{eq:u}) becomes:
\begin{eqnarray}
\label{eq:uu}
U_{\widetilde{m}_0} = Q_{m_0} V_{m_0\times \widetilde{m}_0} \Gamma^{-1/2}_{ \widetilde{m}_0\times \widetilde{m}_0}, 
\hspace*{10mm} \widehat{U}_{\widetilde{m}_0} = \widehat{Q}_{m_0} 
\widehat{V}_{m_0\times \widetilde{m}_0}{\Gamma}^{-H/2}_{ \widetilde{m}_0\times \widetilde{m}_0}.
\end{eqnarray}

\begin{figure}[htbp]
 \center{\includegraphics[width=0.85\textwidth]{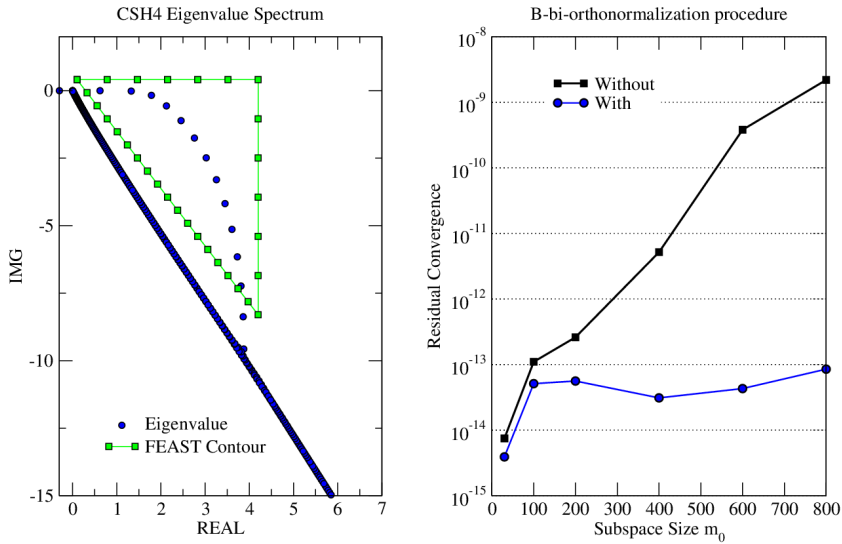}}
        \caption{\label{fig_csh}  On the left: eigenvalue spectrum of CSH4. On the right: minimum 
obtained convergence of the residual norm (\ref{eq:res}) after 20 FEAST iterations plotted
 in function of the subspace size $m_0$. With our bi-orthonormalization procedure, 
the minimum obtained convergence stays relatively constant for all $m_0$.}
\end{figure}

Thereafter, the matrices of the reduced system can be obtained using a new Rayleigh-Ritz projection for $A$ and $B$ i.e. 
 ${B}_U = \widehat{U}^H_{\widetilde{m}_0} B U_{\widetilde{m}_0}$ and ${A}_U = \widehat{U}^H_{\widetilde{m}_0} A U_{\widetilde{m}_0}$.
In spite of our $B$-bi-orthonormalization procedure, the resulting ${B}_U$ is not
 necessarily identity, or even diagonal, due to numerical inaccuracies and finite precision arithmetic.
 However, this procedure is beneficial as a precursor to the QZ algorithm used to solve the reduced 
generalized problem, since 
 it helps to remove contaminating eigenvalues that lie close to 
the contour.
The benefits of our $B$-bi-orthonormalization step can be seen in Figure \ref{fig_csh}. This test has been run on the CSH4 matrix \cite{cerioni2013accurate}, 
an $801 \times 801$ complex scaled Hamiltonian from the BigDFT electronic structure code \cite{genovese2011daubechies}. 
The eigenspectrum and the desired eigenvalues inside of a FEAST custom contour can be seen on the left side of Figure \ref{fig_csh}.
 One edge of the contour is parallel to the eigenvalue continuum. This results in a large number of mixed states 
after spectral projections 
in (\ref{eq:q00}) and (\ref{eq:qhat0}). Without bi-orthonormalization, the QZ algorithm fails to return a $B$-bi-orthogonal 
set of eigenvectors for large values of $m_0$.
 The minimum obtained convergence then degrades for larger subspace sizes. 
By employing our bi-orthogonalization procedure the QZ algorithm is more stable and is able to return a $B$-bi-orthogonal set. 
The minimum obtained convergence remains constant for all $m_0$ values as shown in Figure \ref{fig_csh} (right plot). 
Note that the $B_Q$ matrix remains non-singular for all values of $m_0$ and no resizing operations have then been performed (i.e. $B_U\equiv B_Q$). 

\subsection{Spurious Solutions}\label{sec:spurious}

In certain situation incorrect eigenvalues, so called spurious solutions, appear inside of the FEAST contour. 
These spurious eigenvalues do not converge. It is important to note that the corresponding spurious eigenvectors
 do not need to be explicitly removed from the search subspace to guarantee that the true solutions will
converge along the FEAST iterations. Spurious solutions could then be flagged a posteriori once FEAST has converged.
The spurious problem, however, leads to the practical issue of devising a suitable convergence test.

In FEAST v2.1 for the Hermitian case using Gauss quadrature along a circle contour, 
the true number of eigenvalues $m$ could be obtained by counting
the eigenvalue of $B_Q$ (see  (\ref{eq:bdiag}) using $\widehat{V}=V$) satisfying 
the condition $|\gamma_i|\leq1/4$ \cite{TangPolizzi14,galgon2012counting}
 (i.e. $|\rho_a(\lambda_i)|\leq1/2$ from (\ref{eq:bdiag1})) which guaranteed that $\lambda_i$ is a true 
eigenvalue  within $[\lambda_{min},\lambda_{max}]$.
Since FEAST v3.0 is allowing for custom contour in the complex plane, it is not possible to perform a similar test by
simply analyzing the values $|\gamma_i|$. 
A new strategy has been developed, which can be used to provide increasingly better 
estimate of the number of true eigenvalue solutions in the search subspace at each the FEAST subsequent iteration.

By definition, if a Ritz eigenpair  ($\lambda_i$, $y_i$, $\widehat{y}_i$) obtained after solving the reduced system 
is a genuine solution of the matrix pencil (A,B), 
then the eigenpair ($\rho(\lambda_i)$, $y_i$, $\widehat{y}_i$) is also a solution
for $\rho(B^{-1}A)$ (\ref{eq:spectral2}) and $\rho(AB^{-1})$  (\ref{eq:spectral2hat}). 
In practice, 
one can perform a comparison 
between a direct calculation of  $\rho(\lambda_i)$ where $\lambda_i$ is the Ritz value, and the value 
$\rho(\lambda_i)$ solution of $\rho(B^{-1}A)y_i\simeq\rho(\lambda_i)y_i$ (which is only approximated 
if the Ritz vectors have not yet converged). A suitable choice for the function $\rho$ should
allow these two values for $\rho(\lambda_i)$ to differ significantly if  $\lambda_i$ is spurious, 
with the condition that
$\widehat{y}_i^HB\rho(B^{-1}A)y_i\simeq\rho(\lambda_i)$ can also be easily calculated.
The choice of the approximate spectral projector $\rho_a^2$ (\ref{eq:rational}) satisfy both conditions.
Using (\ref{eq:q00}) and (\ref{eq:bdiag1}), we note that:

\begin{equation}
\label{eq:spurious0}
\rho_a^2(\lambda_i)\simeq\widehat{y}_i^HB\rho_a^2(B^{-1}A)y_i \ = \ \widehat{y}_i^HB{X}\rho_a^2(\Lambda)\widehat{X}^{H}B{y}_{i} \ \equiv \ [B_Q]_{i,i},
\end{equation}
where $[B_Q]_{ii}$ denotes the $i^{th}$ diagonal element of $B_Q$.
Our identification procedure for the spurious solutions can then be summarized by the following three steps:
\begin{enumerate}
\item Compute the corresponding $\{\rho_a(\lambda_i)\}_{i=1,\dots,m_0}$ using (\ref{eq:rational}) and the Ritz values solution of the reduced system $\{\lambda_i\}_{i=1,\dots,m_0}$; i.e.
\begin{eqnarray}
\rho_a(\lambda_i) = \sum_{j=1}^{n_e} \frac{\omega_{j}}{z_j - \lambda_i},
\end{eqnarray}
\item Form the Ritz vectors and wait for the contour integration to be performed and $B_Q$ constructed at the next FEAST iteration.
\item Compare the  calculated values of $\rho_a(\lambda_i)$ with the corresponding diagonal values of $[B_Q]_{ii}$ (which are already sorted),
and label $\lambda_i$ as spurious if it satisfies the following inequality:
\begin{equation}
\left|\frac{\rho_a^2(\lambda_i)-[B]_{ii}}{\rho_a^2(\lambda_i)}\right| \geq \mu,
\end{equation}
where $\mu$ is empirically chosen to be $10^{-1}$. We have found that this criteria is both
large enough to flag all the spurious solutions, and small enough to ensure that 
true solutions are not mislabeled as soon as they start converging.
\end{enumerate}
Once a Ritz eigenpair is flagged as spurious, it is kept in the search subspace but it is not accounted for 
in the test for the residual convergence (\ref{eq:res}). On exit, however, a sorting procedure on the subspace is used by 
FEAST to return the converged eigenpairs free from spurious solutions.

\subsection{Summary and Complete Algorithm}

\begin{figure}[htbp]
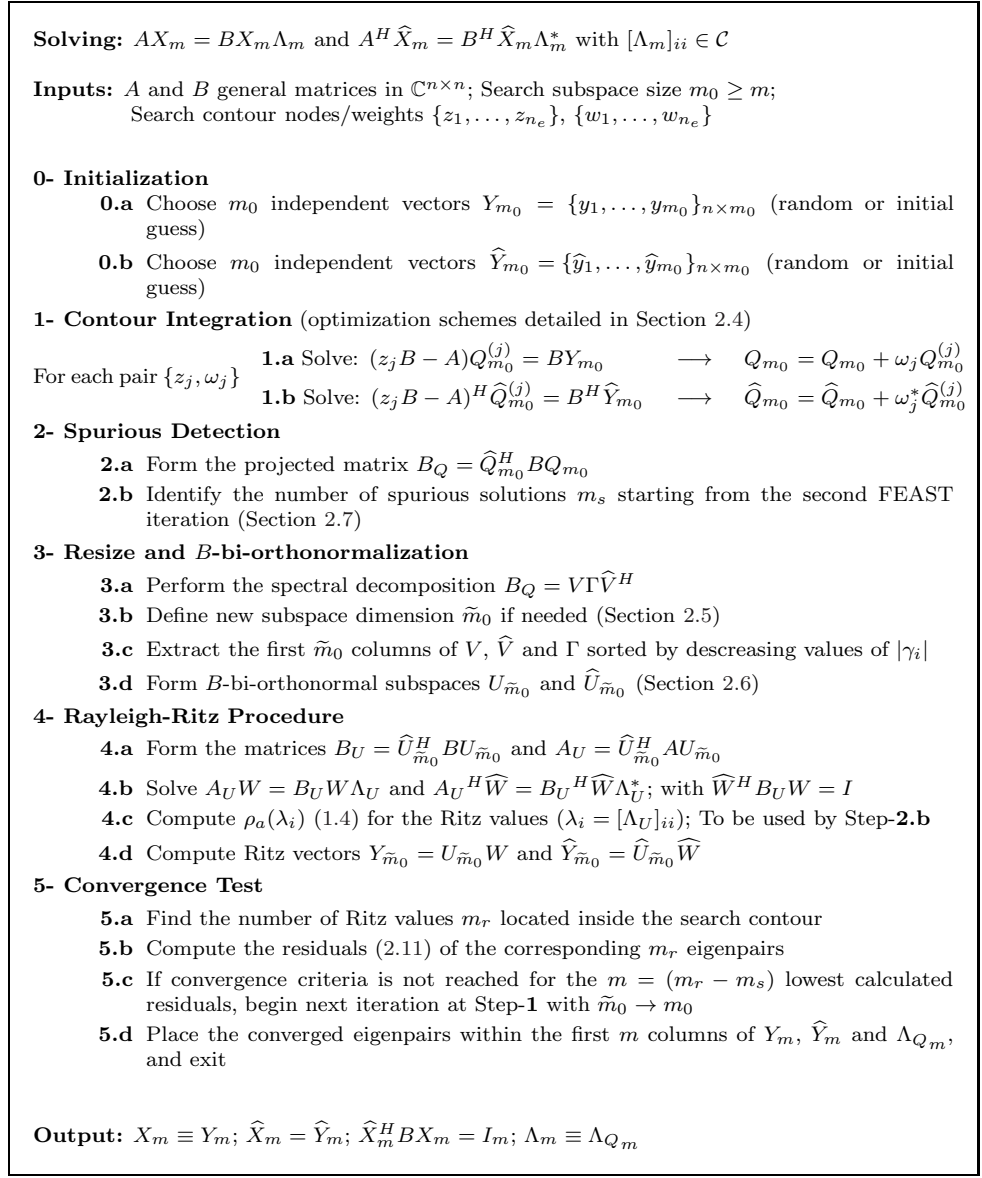

\fbox{
\begin{footnotesize}
\begin{minipage}{0.94\textwidth}
\vspace{2mm}

{\bf Solving:} $AX_m=BX_m\Lambda_m$  and $A^H\widehat{X}_m=B^H\widehat{X}_m\Lambda_m^*$ with $[\Lambda_m]_{ii} \in \cal C$\\

{\bf Inputs:} {$A$ and $B$ general matrices in $\mathbb{C}^{n\times n}$;  Search subspace size $m_0\geq m$; \\ \hspace*{12mm} Search contour nodes/weights 
$\{ z_1,\dots,z_{n_e} \}$, $\{ w_1,\dots,w_{n_e} \}$}
~\\[5pt]

{\bf 0- Initialization}

\begin{itemize}[leftmargin=1.5cm]
\itemsep 1pt
\parskip 1pt
\item [\bf 0.a]  Choose $m_0$ independent vectors $ {Y}_{m_0}=\{y_1,\dots,y_{m_0}\}_{n\times m_0}$ (random or initial guess)
\item [\bf 0.b]  Choose $m_0$ independent vectors ${\widehat{Y}_{m_0}=\{\widehat{y}_1,\dots,\widehat{y}_{m_0}\}_{n\times m_0}}$ (random or initial guess)
\end{itemize}

\vspace{0.1cm}
{\bf 1- Contour Integration} (optimization schemes detailed in Section \ref{sec:contour})
\vspace{0.1cm}

For each pair $\{ z_j,\omega_j\}$ 
	\begin{tabular}{lcl}
	{\bf 1.a} Solve: $(z_j{B} - {A}){Q}^{(j)}_{m_0} = BY_{m_0}$ 
& $\longrightarrow$ & ${Q_{m_0} = Q_{m_0}} + \omega_j {Q}^{(j)}_{m_0}$ \\[3pt] 
	 {\bf 1.b} Solve: $(z_j{B} - {A})^H {\widehat{Q}}^{(j)}_{m_0} = {B^H\widehat{Y}_{m_0}}$  &
$\longrightarrow$ & ${\widehat{Q}_{m_0} = \widehat{Q}_{m_0}} + \omega_j^* {\widehat{Q}}^{(j)}_{m_0}$ 
	\end{tabular}

\vspace{0.1cm}
{\bf 2- Spurious Detection}
\vspace{0.1cm}

\begin{itemize}[leftmargin=1.5cm]
\itemsep 1pt
\parskip 1pt
\item [\bf 2.a] Form the projected matrix $B_Q = \widehat{Q}^H_{m_0} B Q_{m_0}$
\item [\bf 2.b] Identify the number of spurious solutions $m_s$ starting from the second FEAST iteration (Section \ref{sec:spurious})
\end{itemize}

\vspace{0.1cm}
{\bf 3- Resize and $B$-bi-orthonormalization}
\vspace{0.1cm}

\begin{itemize}[leftmargin=1.5cm]
\itemsep 1pt
\parskip 1pt
\item [\bf 3.a] Perform the spectral decomposition $B_Q = V \Gamma \widehat{V}^H$ 
\item [\bf 3.b] Define new subspace dimension $\widetilde{m}_0$ if needed (Section \ref{sec:resize})
\item [\bf 3.c] Extract the first $\widetilde{m}_0$ columns 
of $V$, $\widehat{V}$ and $\Gamma$ sorted by descreasing values of $|\gamma_i|$ 
\item [\bf 3.d] Form $B$-bi-orthonormal subspaces $U_{\widetilde{m}_0}$ and  $\widehat{U}_{\widetilde{m}_0}$ (Section \ref{sec:biortho})
\end{itemize}

\vspace{0.1cm}
{\bf 4- Rayleigh-Ritz Procedure}
\vspace{0.1cm}

\begin{itemize}[leftmargin=1.5cm]
\itemsep 1pt
\parskip 1pt
\item [\bf 4.a] Form the matrices $B_U = \widehat{U}^H_{\widetilde{m}_0} B U_{\widetilde{m}_0}$ 
and $A_U = \widehat{U}^H_{\widetilde{m}_0} A U_{\widetilde{m}_0}$
\item [\bf 4.b] Solve $A_U W =  B_U W \Lambda_U$ and ${A_U}^H\widehat{W}={B_U}^H \widehat{W}\Lambda_U^*$; with $\widehat{W}^HB_UW=I$
\item [\bf 4.c] Compute $\rho_a(\lambda_i)$ (\ref{eq:rational}) for the Ritz values
 ($\lambda_i=[\Lambda_U]_{ii}$); To be used by \mbox{Step-{\bf 2.b}} 
\item [\bf 4.d] Compute Ritz vectors $Y_{\widetilde{m}_0}=U_{\widetilde{m}_0}W$ and $\widehat{Y}_{\widetilde{m}_0}=\widehat{U}_{\widetilde{m}_0}\widehat{W}$
\end{itemize}

\vspace{0.1cm}
{\bf 5- Convergence Test}
\vspace{0.1cm}

\begin{itemize}[leftmargin=1.5cm]
\itemsep 1pt
\parskip 1pt
\item [\bf 5.a] Find the number of Ritz values $m_r$ located inside the search contour
\item [\bf 5.b] Compute the residuals (\ref{eq:res}) of the corresponding $m_r$ eigenpairs  
\item [\bf 5.c] If convergence criteria is not reached for the $m=(m_r-m_s)$ lowest calculated residuals, 
 begin next iteration at Step-{\bf 1} with $\widetilde{m}_0\rightarrow {m}_0$
\item [\bf 5.d] Place the converged eigenpairs within the first $m$ columns of $Y_{{m}}$, $\widehat{Y}_{m}$ and ${\Lambda_Q}_m$,  and exit

\end{itemize}
~\\

{\bf Output:} {$X_m\equiv Y_m$; $\widehat{X}_m=\widehat{Y}_m$; $\widehat{X}_m^HBX_m=I_m$; $\Lambda_m\equiv{\Lambda_Q}_m$}~\\
\end{minipage}
\end{footnotesize}
}
\caption{FEAST Non-Hermitian general algorithm}
\label{alg:three}
\end{figure}

The algorithm in Figure \ref{alg:three} provides a complete description of non-Hermitian FEAST. The algorithm presents
six stages from initialization to convergence test, that further detail the different 
numerical operations outlined in Figure \ref{fig:2algo}. If the non-Hermitian eigenvalue problem is non-defective,
FEAST is expected to converge and return the wanted eigenvalues associated with the $B$-bi-orthonormal right and left eigenvector
subspaces. The convergence rate that was discussed in Section \ref{sec:convergence} depends on the quality of the filter to approximate
spectral projector, and the size of the search subspace (hence it depends on the number of the contour points $n_e$, and 
subspace size $m_0$). Some of the current limitations of the algorithm are outlined in the following:
\begin{description}
\itemsep 1pt
\parskip 1pt
\item[Ill-conditioned linear systems]- In contrast to Hermitian FEAST which allows the selection of  complex shifts (contour points) that 
are not located on the real axis, some of these shifts  could potentially come close to a complex eigenvalue using non-Hermitian FEAST.
Similar to a traditional (Hermitian or non-Hermitian) Arnoldi algorithm using shift-and-invert strategy, 
the resulting linear systems may become ill-conditioned.  
If the shift happens
to be at the exact position of the eigenvalue, the linear system will also be singular. 
One practical solution of this problem consists of moving the contour nodes appropriately
 and automatically by analyzing the eigenspectrum on-the-fly.

\item[Defective system]- Currently if the system is defective, the QZ algorithm used
to solve the reduced system in Step-4b of Figure \ref{alg:three} would not produce a set of $B$-bi-orthogonal subspaces.
In practice, the algorithm may still be found to converge (without Step-2), but 
further studies are required to analyze the action of the approximate spectral projector
on the Jordan form  (\ref{eq:jordan0}) and (\ref{eq:jordan1}).

\item[Ill-conditioned eigenvalue problem]- Non-Hermitian systems are sensitive to the conditioning of the eigenvalues \cite{Golub12}.
A well-known case is the real non-symmetric Grcar matrix \cite{trefethen1991pseudospectra,pseudospectrumgateway} 
(e.g. with $n=100$), which gives rise to extremely sensitive eigenvalues.
It appears some noticeable differences in the eigenvalue calculated using LAPACK-MATLAB, while comparing between 
the eigenvalue solutions of the matrix and its transpose. If double precision arithmetic is desired, this problem 
would require to perform the numerical operations in quad-precision \cite{mct2015}. Interestingly, when FEAST operates on the Grcar matrix or its transpose,
the problem of sensitivity of the eigenvalues is not observed 
in any selected 
regions of the complex plane. 
For this matrix case, the projected reduced eigenvalue problem is then likely to be better conditioned than the original one.
On the other hand, we have found that 
enforcing the condition of bi-orthogonality  could affect the FEAST convergence  for some other systems
e.g. see the case of the QC2534 matrix discussed in Ref. \cite{TangKestynPolizzi14}. 
Further studies are clearly needed to better understand the effect of ill-conditioned eigenvalue systems on FEAST. 

\end{description}

\section{FEAST Eigensolver v3.0 Outlook}\label{sec:3}

The FEAST numerical library package \cite{FEASTsolver} has originally been developed to address
 the Hermitian eigenvalue problem. The package was first released 
(under free BSD license) in Sep. 2009 (v1.0), followed by upgrades in Mar. 2012 (v2.0), and Feb. 2013 (v2.1).
The latter was adopted by Intel math kernel library (Intel-MKL).
 The current version of the FEAST package (v3.0) released in Jun. 2015, started including  
all the various implementation of the non-Hermitian algorithm (real non-symmetric, complex symmetric, and complex general) on both shared-memory systems
(i.e. FEAST-SMP version) and distributed architectures (i.e. FEAST-MPI version).
FEAST's implementation exploit a key strength of modern computer architectures, namely, 
multiple levels of parallelism. FEAST-MPI includes the three levels of
parallelism: MPI for the search contour - MPI for the distribution of the linear
 systems along the contour nodes - OpenMP for the system solver. 

All functionalities of FEAST are accessible through a set of standard predefined interfaces.
The ``ready-to-use'' default drivers are capable to accept dense, banded, and sparse (CSR) matrix formats.
For solving the shifted linear systems, the dense, banded, and sparse FEAST interfaces make use of LAPACK \cite{anderson1999lapack},
 SPIKE-SMP \cite{spike-smp1}, and Pardiso \cite{schenk2006fast} (MKL-version), respectively.
For more advanced users, the FEAST library also includes features
such as reverse communication interfaces (RCI) that are both matrix format and linear system solver independent.
These RCI interfaces 
can then be customized by the end users to
 allow maximum flexibility for their applications. 
In particular, the user is in control of the three major numerical computations to perform on matrices: 
(i) Factorize $(z_jB-A)$ (and $(z_jB-A)^H$ if needed); 
(ii) Solve $(z_jB-A)Q^{(j)}_{m_0}=BY_{m_0}$ and  $(z_jB-A)^H\widehat{Q}^{(j)}_{m_0}=B^H\widehat{Y}_{m_0}$;
(iii)  Mat-vec procedure  involving the multiplications of matrices $A$, $B$, $A^H$, $B^H$ with $m_0$ multiple vectors. 
In order to  address very large sparse systems,
 customized routines such as iterative linear system solvers with or without preconditioners, or
domain decomposition techniques, can straightforwardly be plugged into the RCI loop to perform these operations.
Consequently, the software package has been very well received by the HPC
and application developers, especially in the electronic structure and nanoelectronics communities 
(e.g. \cite{birner2007nextnano,genovese2011daubechies,quantumwise}). 

In addition to the  non-Hermitian interfaces, various supporting routines 
 have also been added in v3.0. These includes: (i) a fast stochastic estimator that can provide a reasonable guess of
the number of eigenvalues count within a user-defined search contour \cite{di2013efficient}; and (ii) a routine 
that can assist the user to extract nodes and weights from a custom design arbitrary geometry in the complex plane.
 This is particularly helpful for non-Hermitian routines as it grants flexibility in targeting specific eigenvalues. 

\section{Numerical Experiments}\label{sec:4}

The non-Hermitian eigenvalue problem (NEP) collection \cite{Bai96} has been used for testing 
and development.
Our test parameters and results for a set of selected system matrices are
 provided  in Table \ref{tab:nep}.  A subset of the eigenpairs has been targeted for each 
system matrix corresponding to
the information provided in the NEP collection, if available. 
Only a few number of FEAST subspace iterations, is needed for most systems to reach convergence.

\begin{table}[htbp]
\begin{center}
\begin{footnotesize}
\begin{tabular}{ l  r r r r  r   c  }
Matrix    &   $n$ & $m_0$ & $m$ & $\lambda_{mid}$ & $r$ & \#Iteration \\
\hline \\ [-1.5ex]
BFW782    & 782 & 44 & 22 & (-5300,300) & 10000.0   & 2 \\
BWM200    & 200 & 36 & 18 & (-1200,0.0) & 60.0 & 2 \\
CDDE5     & 961 & 140 & 70 & (4.75,0.0) & 0.25   & 2\\
GRCAR     & 100 & 38 & 19 & (0.3,0.2) & 0.5  & 4\\
QC324     & 324 & 72 & 37 & (0.0,0.0) & 0.04  & 3 \\
RBS480     & 480 & 112 & 56 & (0.0,0.5) & 0.5   & 9 \\
RW136     & 136 & 38 & 19 & (1.0,0.0) & 0.5  & 5\\
TOLS340   & 340 & 16 & 8 & (-60,300) & 30.0 & 3\\
TOLS4000   & 4000 & 144 & 72 & (-60,300) & 233.0  & 8 \\
\hline
\end{tabular}
\caption{\label{tab:nep}  Non-Hermitian test cases from NEP 
collection. The contour is chosen as a full circle defined by the center and radius $(\lambda_{mid}, r)$
using $n_e=16$ integration points, and the criteria of convergence for the residual is set 
at $10^{-12}$. The system size $n$, the 
subspace size $m_0$, the final number of eigenvalues $m$ found within the
search contour, 
the final residual, and number of FEAST iterations to reach convergence are also listed.}
\end{footnotesize}
\end{center}
\end{table}

\subsection{Parallelism}

As mentioned previously, a major advantage to FEAST are the multiple levels of parallelism naturally 
contained within the algorithm. The following results were gathered on a shared memory machine with 
8 10-core Intel Xeon E7-8870 processors. Each MPI process uses 5 cores.

\begin{figure}[htbp]
\begin{center}
 \center{\includegraphics[width=1.0\textwidth]
        {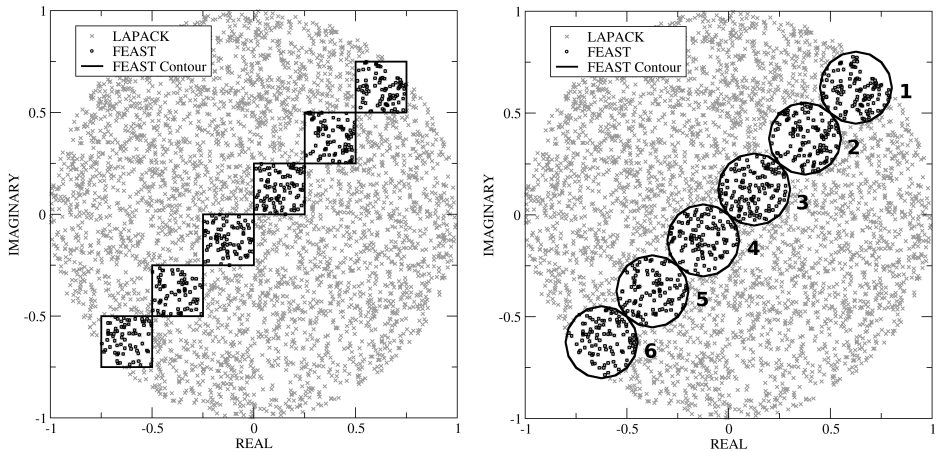}}
\end{center}
\begin{center}
\end{center}
 \caption{\label{fig_contourgrid} A $4000 \times 4000$ dense matrix has been constructed such that all eigenvalues exist within the unit disk. Multiple FEAST contours have been used to calculate a subset of the eigenvalues in parallel. }
\end{figure}

Multiple contours can be solved independently using the 
first level of parallelism of FEAST (overall orthogonality is also largely preserved  \cite{TangPolizzi14,galgon2011feast}). 
However, there is a threshold on the  number of eigenvalues that can be calculated efficiently using a single FEAST contour.
In practice $m_0$ should represent only a small percentage of the matrix size and it may not be suitable to go beyond few thousands
because of the $O(m_0^3)$ complexity of the reduced system solve.
If enough parallel resources are available, however, the solution for an arbitrary large number of eigenvalues can be obtained
by partitioning the entire search domain into multiple
 contours. 
FEAST can then be applied to each  in parallel with a reduced value for $m_0$.
An example of such partitioning is illustrated in Figure \ref{fig_contourgrid}.
The test uses the FEAST dense interfaces on a  $4000 \times 4000$ dense matrix constructed 
such that all eigenvalues exist within the unit disk.
Two sets of contours are considered: 
First, squares with $4$ trapezoidal intervals along each line segment for a total of $16$ linear systems to be solved; 
Next, circles  defined by  16 integration nodes.
In all cases the size of the search subspace is  set at $m_0=200$, and the criteria of convergence
for the residual at $10^{-12}$. At first we consider using only one MPI process per contour, so
the $16$ linear systems are solved one after another using the LAPACK dense solver.
Table \ref{tab_contourgrid} reports the number of eigenvalues found in each contour,  
the number of FEAST iterations, and the total simulation times. Two simulation times are given, 
the fastest has been obtained using a new option offered in FEAST v3.0 that allows to save and reuse 
the factorization at each iteration (increasing then the 
 memory footprint by the number of integration nodes, but removing the need to perform this costly step multiple times).
Saving the factorization between FEAST iterations produced a $2-3\times$ speed improvement for all contours.
As it can be observed from the number of FEAST iterations and the simulation times in Table \ref{tab_contourgrid}, 
load balancing becomes an issue with some contours taking more than twice the time
 of the fastest converging contour.
Since FEAST runs in parallel, its overall efficiency depends on the slowest converging contour (i.e Square 5 or Circle 3). 

\begin{table}[htbp]
\begin{center}
\begin{footnotesize}
\begin{tabular}{| c || c c c c  |}
\hline
Contour $\rm N^o$  & m & \#Iterations & Time-1 (s)  & Time-2 (s)   \\
\hline \hline
 {\bf Square} & &    & &  \\
1        &  84 & 9 & $556$  & $236$ \\
2        & 85  & 7 &$443$  & $197$  \\
3        &  95 & 15 & $891$  & $359$  \\
4        &  83 & 12 & $723$  & $299$  \\
5        & 73  & 19 & $\bf 1107$  & $\bf 438$  \\
6        &  69 & 12 & $718$  & $297$  \\
\hline \hline
 {\bf Circle} & & & &   \\
1        & 120 & 4 &$277$  & $137$   \\
2        & 129 & 8 &$500$  & $217$  \\
3        & 137 & 11 &$\bf 666$ & $\bf 278$    \\
4        & 118 & 8 &$503$   & $218$  \\
5        & 109 & 6 &$389$  & $177$ \\
6        & 104 & 4 &$274$  & $137$  \\ \hline
\end{tabular}
\end{footnotesize}
\end{center}
 \caption{\label{tab_contourgrid} 
 Timing results, number of eigenvalue $m$ and number of iterations obtained using FEAST for each contour 
in Figure \ref{fig_contourgrid}, with $m_0=200$, $n_e=16$ and one single MPI process per contour. 
Two total times are reported by contour: Time-1 for FEAST normal use, and Time-2
that does not account for the cost of the multiple matrix factorizations along the FEAST iterations which are saved in memory.
We note that the overall parallel FEAST efficiency is limited by the slowest individual performance on a single contour
 obtained here for either Square 5 or Circle 3.   }
\end{table}

Better performances can be achieved by taking advantage of another level of parallelism for 
solving the set of independent linear systems.
In the general case, as mentioned in Section \ref{sec:contour}, a single factorization and two solves
 must be performed at each integration node.
With a total of $n_e$ factorizations and $2n_e$ solves, the simulation time could then potentially 
be reduced by a  
factor $n_e$ or more (since the linear systems do not need to be re-factorized at each iteration if $n_e$ is equal to the \#MPI processes).
Table \ref{tab_mpi} presents scalability results for the $4000 \times 4000$ dense matrix considered in Table \ref{tab_contourgrid}.
\begin{table}[htbp]
\begin{center}
\begin{footnotesize}
\noindent
\begin{tabular}{| r ||rrrrrr|| r |} \hline
Contour $\rm N^o$ &     1 & 2 & 3 & 4 & 5 &6 &  Speed-up  \\
\hline
\hline
   {\bf Square}       &          \multicolumn{6}{|c||}{}   &    \\ 
\small 1 MPI & 556 & 443 & 891 & 723 & \bf 1107 &  718  & \bf 1.00 \\
  \small 2 MPI  &  303   & 231   &  457 & 370 &  \bf 566 &  368 &  \bf 1.96    \\ 
  \small 4 MPI   &  160   & 121  &  244 & 198 & \bf 303   & 196  &   \bf 3.65   \\ 
  \small 8 MPI   &  98 & 88 & 149  & 132 & \bf 201 &  128  &     \bf 5.51  \\ 
  \small 16 MPI  &  56  & 41  & 70 & 62 & \bf 95 &69 &   \bf 11.65   \\  \hline \hline
  {\bf Circle}  &     &  & & &  &  &       \\  
\small 1 MPI & 277 & 500 & \bf 666 & 503 & 389 & 274  & \bf 1.00 \\
  \small 2 MPI  &  147 & 252 & \bf 338 & 253 & 196 & 139  &  \bf 1.97  \\ 
  \small 4 MPI  &  81 & 139 & \bf 187 & 140 & 108 & 77 &   \bf 3.56      \\ 
  \small 8 MPI   & 49 & 109  & \bf 148 & 111 & 85 & 60 &  \bf 4.50    \\  
  \small 16 MPI  & 28  & 46 & \bf 60 & 45 & 38 & 30 &    \bf 11.10 \\ \hline
\end{tabular}
\end{footnotesize}
\end{center}
 \caption{\label{tab_mpi} MPI scalability results for the system matrix and contours considered in Figure \ref{fig_contourgrid}
and Table \ref{tab_contourgrid}. The first column indicated the cluster of MPI processes
being used by each contour to distribute the linear systems. 
The last column indicates the speed-up performance associated with the slowest contour (Square 5 or Circle 3).}
\end{table}
For this small  dense example, one observes  only a maximum of 
$\sim11\times$ speed-up compared to a single process  using 16 MPI processes.
The relatively small size  of the test matrix is a limiting factor since
 it leads  to comparable times between
 solving a single linear system and the other numerical operations that 
take place in FEAST (e.g. inner product to form the reduced system,
solution of reduced system, etc.).  
Better scalability performances could be expected using much larger sparse systems.

\section*{Conclusion}
The detailed work developing the non-Hermitian FEAST algorithm has been presented. 
This constitutes a generalization of the well established FEAST Hermitian algorithm,
 leading to a significant upgrade of the FEAST solver package.
The major differences between the Hermitian and non-Hermitian FEAST algorithms stem from the 
complex eigenvalues, which require a two-dimensional search contour. 
Dual subspaces are necessary to allow for computation 
of a $B$-bi-orthogonal basis of left and right eigenvectors. 
In order to improve the stability of the 
algorithm, techniques of subspace resizing, $B$ bi-orthonormalization procedure and spurious detection
 have been implemented and successfully tested.  
We note that the convergence property and parallel capability  associated with the traditional FEAST algorithm have been 
retained with the non-Hermitian algorithm. Finally, the detailed and complete non-Hermitian 
FEAST algorithm implemented in v3.0 is provided, and limitations of its applicability 
 have also been discussed.

\bibliographystyle{siam}
\bibliography{f3}

\begin{thebibliography}{10}

\bibitem{anderson1999lapack}
{\sc E.~Anderson, Z.~Bai, C.~Bischof, S.~Blackford, J.~Demmel, J.~Dongarra,
  J.~Du~Croz, A.~Greenbaum, S.~Hammerling, A.~McKenney, and {and} others}, {\em
  LAPACK Users' guide}, vol.~9, Siam, 1999.

\bibitem{Bai96}
{\sc Z.~Bai, D.~Day, J.~Demmel, and J.~Dongarra}, {\em A test matrix collection
  for non-hermitian eigenvalue problems}, 1996.

\bibitem{bai1999able}
{\sc Z.~Bai, D.~Day, and Q.~Ye}, {\em Able: an adaptive block lanczos method
  for non-hermitian eigenvalue problems}, SIAM Journal on Matrix Analysis and
  Applications, 20 (1999), pp.~1060--1082.

\bibitem{bai1997algorithm}
{\sc Z.~Bai and G.~W. Stewart}, {\em Algorithm 776: Srrit: a fortran subroutine
  to calculate the dominant invariant subspace of a nonsymmetric matrix}, ACM
  Transactions on Mathematical Software, 23 (1997), pp.~494--513.

\bibitem{baker2009anasazi}
{\sc C.~G. Baker, U.~L. Hetmaniuk, R.~B. Lehoucq, and H.~K. Thornquist}, {\em
  Anasazi software for the numerical solution of large-scale eigenvalue
  problems}, ACM Transactions on Mathematical Software (TOMS), 36 (2009),
  p.~13.

\bibitem{birner2007nextnano}
{\sc S.~Birner, T.~Zibold, T.~Andlauer, T.~Kubis, M.~Sabathil, A.~Trellakis,
  and P.~Vogl}, {\em Nextnano: general purpose 3-d simulations}, Electron
  Devices, IEEE Transactions on, 54 (2007), pp.~2137--2142.

\bibitem{cerioni2013accurate}
{\sc A.~Cerioni, Genovese, I.~Duchemin, and T.~Deutsch}, {\em Accurate complex
  scaling of three dimensional numerical potentials}, The Journal of chemical
  physics, 138 (2013), p.~204111.

\bibitem{di2013efficient}
{\sc E.~Di~Napoli, E.~Polizzi, and Y.~Saad}, {\em Efficient estimation of
  eigenvalue counts in an interval}, arXiv preprint arXiv:1308.4275,  (2013).

\bibitem{pseudospectrumgateway}
{\sc M.~Embree and L.~N. Trefethen}, {\em Pseudospectra gateway}.
\newblock \url{http://www.comlab.ox.ac.uk/pseudospectra}.

\bibitem{FEASTsolver}
{\em {\sc FEAST} {\sc e}igenolver}, 2009--2015.
\newblock http://www.feast-solver.org/.

\bibitem{galgon2011feast}
{\sc M.~Galgon, L.~Kr{\"a}mer, and B.~Lang}, {\em The {\sc feast} algorithm for
  large eigenvalue problems}, PAMM, 11 (2011), pp.~747--748.

\bibitem{galgon2012counting}
\leavevmode\vrule height 2pt depth -1.6pt width 23pt, {\em Counting eigenvalues
  and improving the integration in the {\sc feast} algorithm}, Preprint
  BUW-IMACM, 12 (2012), p.~22.

\bibitem{garbow1978algorithm}
{\sc B.~S. Garbow}, {\em Algorithm 535: The qz algorithm to solve the
  generalized eigenvalue problem for complex matrices [f2]}, ACM Transactions
  on Mathematical Software (TOMS), 4 (1978), pp.~404--410.

\bibitem{genovese2011daubechies}
{\sc L.~Genovese, B.~Videau, M.~Ospici, T.~Deutsch, S.~Goedecker, and
  J.~M{\'e}hautois}, {\em Daubechies wavelets for high performance electronic
  structure calculations: The bigdft project}, Comptes Rendus M{\'e}canique,
  339 (2011), pp.~149--164.

\bibitem{Golub12}
{\sc G.~H. Golub and C.~F. Van~Loan}, {\em Matrix computations}, vol.~3, JHU
  Press, 2012.

\bibitem{guettel2014zolotarev}
{\sc S.~G\"uttel, E.~Polizzi, P.~T.~P. Tang, and G.~Viaud}, {\em Zolotarev
  quadrature rules and load balancing for the {\sc feast} eigensolver}, to
  appear in SIAM Journal on Scientific Computing (SISC),  (2015).
\newblock arXiv preprint arXiv:1407.8078 (2014).

\bibitem{hernandez2005survey}
{\sc V.~Hernandez, J.~E. Roman, A.~Tomas, and V.~Vidal}, {\em A survey of
  software for sparse eigenvalue problems}, Universitat Politecnica de
  Valencia, Tech. Rep. STR-6,[retrieved: May, 2013].[Online]. Available:
  http://www. grycap. upv. es/slepc,  (2005).

\bibitem{hernandez2005slepc}
{\sc V.~Hernandez, J.~E. Roman, and V.~Vidal}, {\em {\sc SLEP}c: A scalable and
  flexible toolkit for the solution of eigenvalue problems}, ACM Transactions
  on Mathematical Software (TOMS), 31 (2005), pp.~351--362.

\bibitem{higham2008functions}
{\sc N.~J. Higham}, {\em Functions of matrices: theory and computation}, Siam,
  2008.

\bibitem{knyazev07blopex}
{\sc A.~V. Knyazev, M.~E. Argentati, I.~Lashuk, and E.~E. Ovtchinnikov}, {\em
  Block locally optimal preconditioned eigenvalue xolvers ({\sc blopex}) in
  hypre and petsc}, SIAM Journal on Scientific Computing, 29 (2007),
  pp.~2224--2239.

\bibitem{Laux2012}
{\sc S.~E. Laux}, {\em Solving complex band structure problems with the {\sc
  feast} eigenvalue algorithm}, Physical Review B, 86 (2012), p.~075103.

\bibitem{lehoucq1996evaluation}
{\sc R.~B. Lehoucq and J.~A. Scott}, {\em An evaluation of software for
  computing eigenvalues of sparse nonsymmetric matrices}, Preprint MCS-P547,
  Argonne National Laboratory, 1195 (1996), p.~5.

\bibitem{Arpack98}
{\sc R.~B. Lehoucq, D.~C. Sorensen, and C.~Yang}, {\em {\sc ARPACK} users'
  guide: solution of large-scale eigenvalue problems with implicitly restarted
  Arnoldi methods}, vol.~6, Siam, 1998.

\bibitem{mct2015}
{\sc Advanpix LLC.}, {\em Multiprecision computing toolbox for matlab}.
\newblock \url{http://www.advanpix.com//}.

\bibitem{marek2014elpa}
{\sc A.~Marek, V.~Blum, R.~Johanni, V.~Havu, B.~Lang, T.~Auckenthaler,
  A.~Heinecke, H.~{,} Bungartz, and H.~Lederer}, {\em The elpa library:
  scalable parallel eigenvalue solutions for electronic structure theory and
  computational science}, Journal of Physics: Condensed Matter, 26 (2014),
  p.~213201.

\bibitem{spike-smp1}
{\sc K.~Mendiratta and E.~Polizzi}, {\em A threaded {\sc spike} algorithm for
  solving general banded systems}, Parallel Computing, 37 (2011), pp.~733 --
  741.

\bibitem{moler1973algorithm}
{\sc C.~B. Moler and G.~W. Stewart}, {\em An algorithm for generalized matrix
  eigenvalue problems}, SIAM Journal on Numerical Analysis, 10 (1973),
  pp.~241--256.

\bibitem{parlett1980symmetric}
{\sc B.~A. Parlett}, {\em The symmetric eigenvalue problem}, vol.~7, SIAM,
  1980.

\bibitem{Polizzi09}
{\sc E.~Polizzi}, {\em Density-matrix-based algorithm for solving eigenvalue
  problems}, Physical Review B, 79 (2009), p.~115112.

\bibitem{quantumwise}
{\sc QuantumWise}, {\em Atomistix toolkit version 13.8.1}.
\newblock \url{www.quantumwise.com}.

\bibitem{Saad89}
{\sc Y.~Saad}, {\em Numerical solution of large nonsymmetric eigenvalue
  problems}, Computer Physics Communications, 53 (1989), pp.~71--90.

\bibitem{Saad92}
\leavevmode\vrule height 2pt depth -1.6pt width 23pt, {\em Numerical methods
  for large eigenvalue problems}, vol.~158, SIAM, 1992.

\bibitem{Sakurai2003}
{\sc T.~Sakurai and H.~Sugiura}, {\em A projection method for generalized
  eigenvalue problems using numerical integration}, Journal of computational
  and applied mathematics, 159 (2003), pp.~119--128.

\bibitem{schenk2006fast}
{\sc O.~Schenk and K.~G{\"a}rtner}, {\em On fast factorization pivoting methods
  for sparse symmetric indefinite systems}, Electronic Transactions on
  Numerical Analysis, 23 (2006), pp.~158--179.

\bibitem{sorensen1992implicit}
{\sc D.~C. Sorensen}, {\em Implicit application of polynomial filters in
  ak-step arnoldi method}, Siam journal on matrix analysis and applications, 13
  (1992), pp.~357--385.

\bibitem{stathopoulos2010primme}
{\sc A.~Stathopoulos and J.~R. McCombs}, {\em {\sc PRIMME}: preconditioned
  iterative multimethod eigensolver—methods and software description}, ACM
  Transactions on Mathematical Software (TOMS), 37 (2010), p.~21.

\bibitem{TangKestynPolizzi14}
{\sc P.~T.~P. Tang, J.~Kestyn, and E.~Polizzi}, {\em A new highly parallel
  non-hermitian eigensolver}, in Proceedings of the High Performance Computing
  Symposium, HPC '14, San Diego, CA, USA, 2014, Society for Computer Simulation
  International, pp.~1:1--1:9.

\bibitem{TangPolizzi14}
{\sc P.~T.~P. Tang and E.~Polizzi}, {\em {\sc FEAST} as a subspace iteration
  eigensolver accelerated by approximate spectral projection}, SIAM Journal on
  Matrix Analysis and Applications, 35 (2014), pp.~354--390.

\bibitem{trefethen1991pseudospectra}
{\sc L.~N. Trefethen}, {\em Pseudospectra of matrices}, Oxford University,
  Computing Laboratory Numerical Analysis Group, 1991.

\bibitem{trefethen2014exponentially}
{\sc L.~N. Trefethen and Weideman J.~A. C.}, {\em The exponentially convergent
  trapezoidal rule}, SIAM Review, 56 (2014), pp.~385--458.

\bibitem{yin2014feast}
{\sc G.~Yin, R.~H. Chan, and M.~Yeung}, {\em A {\sc feast} algorithm for
  generalized non-hermitian eigenvalue problems}, arXiv preprint
  arXiv:1404.1768,  (2014).

\bibitem{zolotarev1877application}
{\sc E.~I. Zolotarev}, {\em Application of elliptic functions to questions of
  functions deviating least and most from zero}, Zap. Imp. Akad. Nauk. St.
  Petersburg, 30 (1877), pp.~1--59.

\end{thebibliography}
%
%
%

\end{document}